\documentclass{article}
\usepackage{amssymb}
\usepackage{amsmath}
\usepackage{amsthm} 
\usepackage{graphicx}
\usepackage{soul} 
\usepackage{color, colortbl} 
\usepackage{subcaption}
\usepackage{float} 

\usepackage{pst-plot}
\usepackage{tikz}

\newtheorem{theorem}{Theorem}
\newtheorem{definition}{Definition}
\newtheorem{lemma}[theorem]{Lemma}

\newtheorem{con}[theorem]{Conjecture}

\newtheorem{example}{Example}

\newcommand{\rank }{\,{\rm rank}}

\newcommand{\block}[2]{
  \underbrace{
  \begin{matrix}
  #2
  \end{matrix}}_{#1}
}

\def\reals{\mathbb{R}}

\def\qnm{Q_n^m}

\def\cong{\rm cong}
\def\iso{\rm iso}

\def\RS{${\cal RS}$}

\newcommand{\Rred}{R^{\mathrm{red}}}
\newcommand{\Fred}{F^{\mathrm{red}}}

\newcommand{\comment}[1]{\mbox{}}

\def\qed{{\hfill{\vrule height5pt width3pt depth0pt}\medskip}}

\begin{document}

\begin{center} {\bf \LARGE Central Configurations with Unequal Masses: Finiteness in Several Exceptional Cases of Five Bodies}
\vskip 0.5cm
Ma{\l}gorzata Moczurad and Piotr Zgliczy\'nski \\
   \{malgorzata.moczurad, piotr.zgliczynski\}@uj.edu.pl \\
Faculty of Mathematics and Computer Science, Jagiellonian University,\\
ul. {\L}ojasiewicza 6,
30-348 Krak\'ow, Poland

\vskip 0.5cm
 \today
\end{center}

\begin{abstract}
We provide a computer-assisted proof of the exact count of classes of central configurations for five bodies for several sets of mass values that are exceptional from the point of view of the finiteness results of Albouy and Kaloshin in the planar case and of Hampton and Jensen in the spatial case.
\end{abstract}

\tableofcontents

\section{Introduction}

A long-standing question, raised by Wintner~\cite{Wintner}, concerns the
finiteness  of the number of  central
configurations for any choice of positive. In~\cite{SmNext}, Smale listed this problem among a number of
challenging problems for the twenty-first century, where it appears as
Problem~6.

The relative equilibria of the three-body problem have been known since the
eighteenth century. Up to equivalence, there are exactly five such
configurations. Three of them are collinear configurations discovered by
Euler, while the remaining two correspond to Lagrange’s equilateral
triangles. Euler’s collinear configurations were later generalized to the
$n$-body problem by Moulton~\cite{Mou}, who showed that there are exactly
$n!/2$ collinear equivalence classes.

The most successful approaches to the finiteness problem for central
configurations exploit the fact that the defining equations can be written as
systems of polynomial equations. Techniques from algebraic geometry and tropical algebraic geometry are then applied to study these systems.
The resulting proofs are computer-assisted, relying on symbolic
computations and/or exact integer arithmetic. Below we summarize the major
developments achieved using this approach.

In 2006, Hampton and Moeckel~\cite{HM} proved that the number of relative
equilibria in the Newtonian four-body problem is finite, lying between 32 and
8472. Their proof is computer-assisted and based on symbolic and exact integer
computations. The upper bound of 8472 is believed to be a significant
overestimate; so far, no more than 50 equilibria have been found (see, for
example,~\cite{Si78}).

In 2012, Albouy and Kaloshin~\cite{AK} nearly resolved the finiteness question
for the planar five-body problem. They proved that there are finitely many relative equilibria, except possibly
when the mass vector belongs to a certain codimension--2 subvariety
$\mathcal{S}$ of the mass space. The key idea of their proof is to follow a
potential continuum of central configurations into the complex domain and to
analyze its possible singularities. An upper bound on the number of relative
equilibria follows from B\'ezout’s theorem, although the authors remark that
\emph{the bound is so bad that we avoid writing it explicitly}. It is worth
noting that the equal-mass case is exceptional in the sense of Albouy and
Kaloshin, that is, it belongs to the set $\mathcal{S}$.

The spatial five-body problem was studied by Hampton and Jensen~\cite{HJ}, who
combined polyhedral and polynomial computations to derive equations describing
the set of exceptional mass choices for which finiteness of central
configurations may fail. Their work generalizes an earlier generic finiteness
result of Moeckel~\cite{M01}.

More recently, Jensen and Leykin~\cite{JL}, as well as Chang and
Chen~\cite{CC24,CC25}, investigated the planar six-body problem. Jensen and
Leykin employed techniques from tropical geometry, while Chang and Chen
implemented the approach of Albouy and Kaloshin in an algorithmic framework.
Both works re-established the generic finiteness result for the planar
five-body problem proved in~\cite{AK}. Attempts to extend these methods to the
case $n=6$ have so far been unsuccessful.

Our paper follows a different line of attack on the finiteness problem, based
on techniques from interval arithmetic~\cite{Mo}. Although this approach has
not yet produced results as far-reaching as those described above, it has led
to proofs of finiteness and complete classification of central configurations
for equal masses in the planar case for $n=4,5,6$ and $7$~\cite{MZ}, and in the
spatial case for $n=4,5$ and $6$~\cite{MZ20}. At first glance, the equal-mass
case might appear to be highly degenerate, suggesting that this approach could
fail for generic mass distributions. However, due to the nature of our
arguments, the resulting counts of distinct classes of central configurations
remain valid for masses lying in a small neighborhood of the equal-mass case.  Let us mention here that the most comprehensive and possibly complete
list of planar central configurations for equal masses for $n$ up to $12$ has been  published in \cite{DZD}.   

The main limitation of these results is that they apply only to specific mass
choices (or small neighborhoods thereof). In this context, it is worth
mentioning the work~\cite{FTZ}, which provides a computer-assisted proof of
finiteness and a complete classification of relative equilibria for the planar
restricted four-body problem when the massive bodies form an equilateral
triangle. In that setting, all positive masses are allowed. From the
perspective of the general finiteness problem, this may be viewed as a toy
model. Nevertheless, it captures several key difficulties that must be
addressed by interval arithmetic methods, including high-dimensional parameter
spaces, bifurcations, and singular limits arising when some masses tend to
zero. We note in passing that the same result was obtained analytically by
Barros and Leandro~\cite{BL11,BL14}.

The present paper is a sequel to~\cite{MZ,MZ20}. We apply our approach to the
five-body problem for several mass distributions that are exceptional from the
point of view of the finiteness results of Albouy and Kaloshin in the planar
case~\cite{AK}, and of Hampton and Jensen in the spatial case~\cite{HJ}. For
these mass choices, we show that all central configurations are non-degenerate
and provide an exact count of them. It is worth emphasizing that the
equal-mass case itself is also exceptional in the sense of~\cite{AK,HJ}. This
suggests that the exceptional character of these cases is more likely due to
limitations of the techniques used in~\cite{AK,HJ}, rather than to intrinsic
properties of the corresponding central configurations.

Although the approach developed in~\cite{MZ} for the planar case and later adapted in~\cite{MZ20} to the spatial case is, in principle, applicable to the case of unequal masses, additional mathematical insights were required to make the method effective. In~\cite{MZ,MZ20}, we employed a subset of the equations defining central configurations—referred to as the reduced system—in which several equations were omitted. This reduced system was fixed and used consistently throughout the computations. In the equal-mass case, symmetry arguments allowed us to restrict the search space for possible central configurations, so that all such configurations turned out to be non-degenerate solutions of the reduced system and were therefore amenable to computer-assisted proofs based on interval arithmetic versions of Newton’s method.

In the case of unequal masses, however, a non-degenerate central configuration may correspond to a degenerate solution of one reduced system, while remaining a non-degenerate solution of another. The main objective of this work is to understand when and how this phenomenon occurs and to describe an implementation that avoids this issue by switching between different reduced systems. To achieve this, we identify the precise conditions under which a non-degenerate central configuration becomes a degenerate solution of a reduced system. In the planar case, this characterization is given in Theorems~~\ref{thm:RS-nondeg-equiv} and~\ref{thm:sol-RS-deg} in Section~\ref{subsec:ndeg2d}, while in the spatial case it is provided in Theorems~\ref{thm:sol-RS3D-deg} and~\ref{thm:RS3D-nondeg-equiv} in Section~\ref{subsec:ndeg3d}.

\section{Reduced system of equations for normalized central configurations}
\label{sec:cc-eq}

Assume there is a group  of $n$ bodies (point masses)  interacting with each other gravitationally (i.e.\ due to the  Newton's law of gravitation; the gravitational constant is normalized $G = 1$).

\begin{definition} \label{def:conf}\rm
Let $q = (q_1, \ldots, q_n)\in(\reals^d)^n$ and $m = (m_1, \ldots, m_n)\in \reals_+^n$.
An element $\qnm = (q, m)$
will be called  a \emph{configuration} (of $n$ bodies). Each body has a position $q_i\in\mathbb{R}^{d}$ and a mass $m_i\in\reals_+$.
The coordinates $q_i$ will be given  as $q_i=(x_i,y_i)$ when $d=2$ or $q_i=(x_i,y_i,z_i)$ when $d=3$.
\end{definition}

The center of mass for configuration $\qnm$ is given by
\begin{equation}
   c=\frac{1}{\sum_i m_i} \sum_{j} m_j q_j.
\end{equation}

The {\em central configuration problem}: for given masses $m_i$  to find   $c \in \reals^d$, $\lambda \in \reals_+$  and positions of bodies  satisfying the following system of equations
\begin{equation}\label{eq:cc-with-lambda}
  \lambda m_i(q_i-c) = \sum_{\substack{j=1\\
j\neq i}}^n \frac{m_im_j}{r_{ij}^3}(q_i - q_j)=: f_i(q_1,\dots,q_n), \quad i=1,\dots,n.
\end{equation}
It turns out that $c$ must be a center of mass of configuration $\qnm$ and $\lambda \in \reals_+$ is also determined by $\qnm$.

 Clearly, no solution of~(\ref{eq:cc-with-lambda}) is isolated, because close to any configuration is another, obtained from the former by rotation, scaling or translation. For this reason, equivalence classes of central configurations are introduced. Two configurations are called equivalent if they can be transformed into each other by
rotating around the center of mass, scaling or shifting accordingly.

The goal of this section  is to present a set of equations ({\em the reduced system of equations}), which gives all equivalence classes of CCs; this is  based on Section~2 in \cite{MZ22}, Section~5 in \cite{MZ} for the planar case and Section~3 in \cite{MZ20}
for the spatial case. First,  we eliminate scaling and translational symmetries simply by setting $\lambda = 1$ and $c = 0$ (Definition~\ref{def:ncc}). Afterwards, we remove an equation for the last body using the center of mass reduction (Section~\ref{subsubsec:com-red}) and finally we
 remove $SO(d)$ symmetry by demanding that selected body is on OX-axis (if $d=2$) (see Section~\ref{sec:red-sys2D}) and when $d=3$  we additionally demand that some other body is on $OXY$ plane (see Section~\ref{sec:red-sys-equiv}).

\begin{definition}\label{def:ncc}\rm\cite{AK,MZ}
A \emph{normalized central configuration} is a solution of (\ref{eq:cc-with-lambda}) with $\lambda=1$ and $c=0$.
\end{definition}

Henceforth, nCC denotes a normalized central configuration, while CC denotes a central configuration.

The system of equations for normalized central configurations is
\begin{equation}\label{eq:cc-kart}
  q_i= \sum_{j,j\neq i} \frac{m_j}{r_{ij}^3}(q_i - q_j)= \frac{1}{m_i}f_i(q_1,\dots,q_n), \quad i=1,\dots,n.
\end{equation}
From now on we focus on normalized central configurations.
We introduce the function $F\colon (\reals^{d })^n \to (\reals^{d })^n$ given by
\begin{equation}\label{eq:vector-field}
  F_i(q_1,\dots,q_n) =  q_i - \sum_{j,j\neq i} \frac{m_j}{r_{ij}^3}(q_i - q_j), \quad i=1,\dots,n.
  \index{$F_i$}
\end{equation}
With this notation~(\ref{eq:cc-kart})
becomes
\begin{equation}\label{eq:cc-abstract}
  F(q_1,\dots,q_n)=0.
  \index{$F$}
\end{equation}

It is well known (see \cite{MZ} and the literature given there) that for any $(q_1,\ldots,q_n)\in (\reals^d)^n$ holds
\begin{eqnarray}
\sum_{i=1}^n f_i&=&0, \label{eq:n-mom-con} \\
\sum_{i=1}^n f_i \wedge q_i & = & 0, \label{eq:n-angular-mom-con}
\end{eqnarray}
where $v \wedge w$\index{$v \wedge w$} is the exterior product of vectors, the result being an element of exterior algebra. If $d=2$ or 3 it can be interpreted as the vector product of $v$ and $w$ in dimension $3$.  The identities (\ref{eq:n-mom-con}) and (\ref{eq:n-angular-mom-con})
are easy consequences of the third newton's law (the action equals reaction) and the requirement that the mutual forces between bodies are in direction of the other body.

\subsection{Center of mass reduction}\label{subsubsec:com-red}

Consider system (\ref{eq:cc-abstract}). After multiplication of  $i$-th equation by $m_i$ and addition of all equations using (\ref{eq:n-mom-con}) we obtain (or rather recover)
the center of mass equation
\begin{eqnarray}
 \left(\sum_{i=1}^n m_i\right) c=\sum_i m_i q_i = 0. \label{eq:cc-cofmass-l}
\end{eqnarray}
We can take the equations for $n$-th body and replace it with (\ref{eq:cc-cofmass-l}) to obtain an equivalent system
\begin{subequations}
\begin{align}
  q_i&= \sum_{j=1,j\neq i}^n \frac{m_j}{r_{ij}^3}(q_i - q_j), \quad i=1,\dots,n-1, \label{eq:cc-kart-1n-1} \\
   q_n&=-\frac{1}{m_n}\sum_{i=1}^{n-1} m_i q_i. \label{eq:cc-kart-n-th} 
\end{align}\label{eqn:reduced}
\end{subequations}

We  write the  system (\ref{eqn:reduced}) obtained from (\ref{eq:cc-abstract}) after removing the $n$-th body using the center of mass equation (condition (\ref{eq:cc-cofmass-l})) as
\begin{equation}
  \Fred(q_1,\dots,q_{n-1})=0,
\end{equation}
where $\Fred: (\reals^d)^{n-1} \to (\reals^d)^{n-1}$.
To be precise we have
\begin{eqnarray*}
   \Fred_i(q_1,\dots,q_{n-1})&=& F_i(q_1,\dots,q_{n-1},q_n(q_1,\dots,q_{n-1})),\quad i=1,\dots,n-1,
\end{eqnarray*}
where
\begin{equation}
    q_n(q_1,\dots,q_{n-1})=-\frac{1}{m_n}\sum_{i=1}^{n-1} m_i q_i. \label{eq:cc-com}
\end{equation}

Let us introduce a  notation
$$R\colon (\reals^d)^n\to (\reals^d)^n\quad \mbox{and}\quad \Rred\colon (\reals^d)^{n-1}\to (\reals^d)^{n-1}$$
that will facilitate the manipulation of the system of equations.
For any configuration $\qnm$ we set
\begin{equation*}
  R_i(q_1,\dots,q_n)=m_iF_i(q_1,\dots,q_n) = m_iq_i - f_i(q_1,\dots,q_n), \quad i=1,\dots,n.
\end{equation*}
With the above notation the system (\ref{eq:cc-abstract}) becomes
\begin{equation*}
  R(q_1,\dots,q_n)=
  (R_1(q_1,\dots,q_n),\dots,R_n(q_1,\dots,q_n))=0. \label{eq:cc-R}
\end{equation*}
For any $(q_1,\dots,q_{n-1}) \in (\mathbb{R}^d)^{n-1}$
we define
\begin{equation*}
  \Rred_i(q_1,\dots,q_{n-1})=R_i(q_1,\dots,q_{n-1},q_n(q_1,\dots,q_{n-1})), \quad i=1,\dots ,n-1.  \label{eq:R-tylda}
\end{equation*}

With the above notation we have
\begin{equation}
  m_i \Fred_i(q)=\Rred_i(q), \quad i=1,\dots,n-1.
\end{equation}
In our previous papers \cite{MZ20,MZ22} we used $\tilde{R}$ for $\Rred$.

The following lemma was proved in \cite{MZ22}.
\begin{lemma}\cite[Lemma 2]{MZ22}
\label{lem:c-angmom-red}\rm
For any $(q_1,\dots,q_{n-1})\in (\mathbb{R}^d)^{n-1}$  holds
\begin{equation*}
  \sum_{i=1}^{n-1}(q_i-q_n(q_1,\dots,q_{n-1})) \wedge \Rred_i(q_1,\dots,q_{n-1}) = 0. \label{eq:c-angmom-red}
\end{equation*}
\end{lemma}


\subsection{Reduced systems}
 For $n$ bodies, we started with a system~(\ref{eq:cc-with-lambda}) of $dn$ equations with  of $dn+2$ unknowns; by fixing the center of mass $c = 0$ and  $\lambda = 1$, we got a system of $dn$ equations with $dn$ unknowns. We further eliminated the equations for one body using the center of mass reduction. So at this stage, for $n$ bodies, we already have only $d(n-1)$ equations with $d(n-1)$ unknowns, but still the system has a rotational symmetry, so its solutions are not isolated.
 In what follows, we derive a reduced system of equations such that each equivalence class of nCCs has exactly one solution. To achieve this, we remove the  $\mathcal{SO}(d)$ symmetry.

\subsubsection{General remarks on the systems of equations introduced in subsequent sections}
In this work, we use two notations for the system of equations defining central configurations:
$F$ and $R$. These notations (both in their full and reduced forms) are interdependent and equivalent, since
$m_iF_i(q) = R_i(q)$. The $R$ notation is used for theoretical results, as it is easier to manipulate. However, in the program used to find central configurations, we solve the equations in terms of $F$; therefore, both notations are mentioned in this work.

Later on we will also introduce a reduced system denoted \RS.
These  systems  differ in the number of variables, i.e.
\begin{eqnarray*}
\Fred, \Rred \colon  \left(\reals^d \right)^{n-1} \to \left(\reals^d \right)^{n-1} & & \\
\mathcal{RS}  \colon \left(\reals^2 \right)^{n-2}\times \reals \to \left(\reals^2 \right)^{n-2}\times\reals & & \qquad \mbox{(in the case 2D  see Definition~\ref{def:RS})} \\
\mathcal{RS}  \colon \left(\reals^3 \right)^{n-3}\times \reals^2\times \reals\to \left(\reals^3 \right)^{n-3}\times\reals^2\times\reals & & \qquad \mbox{(in the case 3D  see Definition~\ref{def:sol-RS-3D})}
\end{eqnarray*}
It is reflected in the types and forms of their Jacobian matrices. For example, in 3D case, $D\Fred(q)$ and $D\Rred(q)$ have $3n - 3$ rows (and columns) whereas $D\mathcal{RS}(q)$ has only $3n - 6$. This affects the number of active variables in a given system of equations; nevertheless, throughout this work, unless this could cause ambiguity, we denote them simply by
$q$, i.e.\, $\mbox{\RS}(q)$ instead of $\mbox{\RS}(q_1, \ldots, q_{n-1})$. We assume that the reader will apply the appropriate type of function in the relevant context.

We introduce the reduced systems \RS\  by eliminating certain equations. The choice of which equations to eliminate is made solely for notational convenience; since bodies may be permuted and configurations rotated, the index of the eliminated equation is irrelevant.

\subsubsection{Reduced system in 2D}\label{sec:red-sys2D}

The goal of this section is to define a reduced system on the plane (i.e.\, $d = 2$).  We follow Section 2.2 in \cite{MZ22}.
We  use the notation  $\Fred_i=(\Fred_{i,x},\Fred_{i,y})$ and $\Rred_i=(\Rred_{i,x},\Rred_{i,y})$.
Let us fix $k_0\in \{1,\dots,n-1\}$, and consider the following set of equations (compare (\ref{eq:cc-kart}))
\begin{subequations}
\begin{align}
  q_i &= \frac{1}{m_i}f_i(q_1,\dots,q_n(q_1,\dots,q_{n-1})), \quad i\in \{1,\dots,n-1\}, i \neq k_0,  \\ 
  x_{k_0} &=  \frac{1}{m_{k_0}}f_{k_0,x}(q_1,\dots,q_n(q_1,\dots,q_{n-1})),
\end{align}\label{eqn:cc-ref}
\end{subequations}
where $f_i = (f_{i,x}, f_{i, y})$.\index{$f_{k,x}$}
Observe that system (\ref{eqn:cc-ref})  has $2(n-1)-1$ equations for $q_1,\dots,q_{n-1} \in \mathbb{R}^2$ and it coincides with~(\ref{eqn:reduced}) with the
equation for $y_{k_0}$ dropped. To obtain the same number
of independent variables we set $y_{k_0}=0$.

Using the notation introduced in Section~\ref{subsubsec:com-red} system (\ref{eqn:cc-ref}) can  be  equivalently written as
\begin{subequations}
\begin{align}
 \Fred_i(q_1,\dots,q_{n-1}) &=0, \quad i\in \{1,\dots,n-1\},\ i \neq k_0, \label{eq:cc-R-red-i} \\
  \Fred_{k_0,x}(q_1,\dots,q_{n-1}) &=0, \label{eq:cc-R-red-xk1}
\end{align}\label{eq:cc-R-red-xy}
\end{subequations}\\
 where we substitute $0$ for $y_{k_0}$.
The next theorem addresses the question: whether from a solution of (\ref{eq:cc-R-red-xy}) we
obtain a solution of (\ref{eq:cc-kart})?

\begin{theorem}\cite[Theorem 3]{MZ22}
\label{thm:red-to-full2D}
   If $\overline{q} = (q_1,\dots,q_{n-1})$ satisfies equations  (\ref{eq:cc-R-red-xy})  and is  such that
    \begin{equation}\label{eqn:A1}
 x_{k_0} \neq x_n,
    \end{equation}
     then
    $q = (\overline{q}, q_n(\overline{q}))$ is a normalized central configuration, i.e.\ it satisfies (\ref{eq:cc-kart}).
\end{theorem}


\begin{definition}\cite[Def. 4]{MZ22}\label{def:RS}\rm
~System of equations  (\ref{eq:cc-R-red-xy}) with $k_0=n-1$ and with $y_{n-1}=0$ will be called {\em the reduced system}. We will use abbreviation \RS\ for this system.
\end{definition}

The system \RS\ no longer has $\mathcal{O}(2)$ as its symmetry group; however, it remains symmetric with respect to reflections across the coordinate axes $OX$ and $OY$.

In defining \RS, we made two arbitrary choices: which body is determined by the center-of-mass condition~\eqref{eq:cc-com} (the $n$-th body), and which body is placed on the $OX$ axis (the $(n-1)$-st body). Both choices are inessential, since permutations of the bodies and rotations of the coordinate system allow any body to play either role.

The variables of \RS\ are $(q_1,\dots,q_{n-2},x_{n-1}) \in \left(\mathbb{R}^2\right)^{n-2} \times \mathbb{R}$, hence its solutions are not configurations of $n$ bodies in the sense of Definition~\ref{def:conf}. However, to facilitate further discussion we introduce the following convention.

\begin{definition}\label{def:sol-RS}\rm
We say that $\qnm$ satisfies \RS (or, informally, is a solution of \RS) iff
$$
\left\{
\begin{array}{rcl}
 \Fred_i(q_1,\dots,q_{n-1}) & = & 0, \quad i\in \{1,\dots,n-1\},\ i \neq n-1, \\[1ex]
  \Fred_{n-1,x}(q_1,\dots,q_{n-1}) & = & 0, \\[1ex]
  y_{n-1} & = & 0,\\[1ex]
  q_n & = & q_n(q_1,\dots,q_{n-1}).
\end{array}
\right.
$$
\end{definition}

We can also define another reduced system by requesting that $y_{n-1}=y_{n-2}$. In such situation  the variable set is the same as above
i.e. $(q_1,\dots,q_{n-2},x_{n-1}) \in \left(\mathbb{R}^2\right)^{n-2} \times \mathbb{R}$, but this time the full configuration  is defined by setting $y_{n-1}=y_{n-2}$. This kind of normalization was used in \cite{AK}.

\begin{theorem}
\label{thm:non-deg-RS}
Let $\overline{q} = (q_1, \ldots, q_{n-1})$ and $q = (\overline{q}, q_n(\overline{q}))$.
 Assume that  $\qnm = (q, m)$   is a  nCC. Then, in a suitable coordinate system and after  some permutation of bodies,  $\overline{q}$ is a solution of \RS\ satisfying $x_n \neq x_{n-1}$ and $x_{n-1} \neq 0$.
\end{theorem}
\proof
 First we take any $q_{i_0}\neq 0$ (there can be only one body at the origin) and we chose coordinate frame so that $q_{i_0}=(x_{i_0},0)$. Then we look for $j \neq i_0$ such that $x_{j} \neq x_{i_0}$.  Observe that due to the center of mass condition (\ref{eq:cc-cofmass-l}) such $j$ always exists.  Now we change the numeration of bodies so that $i_0 \to n-1$ and $j \to n$.
\qed

In Theorem~\ref{thm:non-deg-RS} the presence of first condition  $x_n \neq x_{n-1}$ is motivated by Theorem~\ref{thm:red-to-full2D}. The second  condition $x_{n-1} \neq 0$ is related to the non-degeneracy question of nCC, see Theorems~\ref{thm:RS-nondeg-equiv} and~\ref{thm:sol-RS-deg} from Section~\ref{subsec:ndeg2d}.

\subsubsection*{Solution of \RS\  which is not an nCC}

Note that being a solution of \RS\ is not sufficient to be an nCC. In \RS\ we omit the equation for $y_{n-1}$ assuming $y_{n-1} = 0$, but in fact this equation has to be satisfied.
\begin{example}\rm
\label{ex:ex1}
Consider a collinear configuration of three bodies lying on the OY-axis (i.e. $x_i=0$ for $i=1,2,3$) and $(y_1>0,y_2=0,y_3=-\frac{m_1}{m_3}y_1)$. Since manifestly, $\Rred_{1,x}=0$ and $\Rred_{2,x}=0$, it will be a solution of \RS\ if  $\Rred_{1,y}=0$ which is equivalent to
 \begin{eqnarray*}
  y_1 - \frac{m_2}{y_1^2} - \frac{m_3}{y_1^2 \left(1+ \frac{m_1}{m_3} \right)^2} & = &0
\end{eqnarray*}
and finally we obtain
\begin{eqnarray}
  m_2 + \frac{m_3}{\left(1+ \frac{m_1}{m_3} \right)^2} & = & y_1^3. \label{eq:rsys-y1}
\end{eqnarray}

Clearly, such a configuration does not satisfy condition~(\ref{eqn:A1}).
We are therefore led to the question of whether, under these assumptions, the equation $\Rred_{y,2}=0$ given by
\begin{eqnarray*}
  y_2 - \frac{m_1 (y_2- y_1)}{|y_2-y_1|^3} - \frac{m_3 (y_2-y_3)}{|y_2-y_3|^3} & = & 0,
  \end{eqnarray*}
 is satisfied. It turns out that this is the case only when $m_1 = m_3$.

 Moreover, it could be shown that if $m_1 \neq m_3$, then it is a non-degenerate solution of \RS\ .

\end{example}
\subsubsection{Reduced system  in 3D}
\label{sec:red-sys-equiv}

The aim of this section is to derive the reduced system of equations in the spatial case, i.e., for $d = 3$. We follow Section 3 in \cite{MZ20}.
Let us fix $k_1,k_2 \in \{1,\dots,n-1\}$, $k_1 \neq k_2$ and consider the following set of equations
\begin{subequations}
\begin{eqnarray}
  q_i &=&\frac{1}{m_i}f_i(q_1,\dots,q_n(q_1,\dots,q_{n-1})), \quad i\in \{1,\dots,n-1\}, i \neq k_1,k_2 \label{eq:cc-red-i} \\
  x_{k_1} &=& \frac{1}{m_{k_1}}f_{k_1,x}(q_1,\dots,q_n(q_1,\dots,q_{n-1})), \label{eq:cc-red-xk1}\\
   x_{k_2} &=& \frac{1}{m_{k_2}}f_{k_2,x}(q_1,\dots,q_n(q_1,\dots,q_{n-1})), \label{eq:cc-red-xk2}\\
    y_{k_2} &=& \frac{1}{m_{k_2}}f_{k_2,y}(q_1,\dots,q_n(q_1,\dots,q_{n-1})), \label{eq:cc-red-yk2}
\end{eqnarray}\label{eq:cc-red-3D}
\end{subequations}
where $f_i = (f_{i,x}, f_{i, y},f_{i,z})$\index{$f_{k,x}$} and
\begin{eqnarray}
  q_n(q_1,\dots,q_{n-1})&=&-\frac{1}{m_n}\sum_{i=1}^{n-1} m_i q_i. \label{eq:cc-red-com}
\end{eqnarray}
As in the planar case ($d = 2$), we denote by \RS\  the reduced system of equations. In the spatial case considered here, \RS\ consists of equations~(\ref{eq:cc-red-3D}), under the assumption  $y_{k_1} = z_{k_1} = z_{k_2} = 0$.

\begin{definition}\rm
System of equations~(\ref{eq:cc-red-3D}) with $k_1=n-1$, $k_2 = n-2$ and with $y_{n-1}= z_{n-1} = z_{n-2} = 0$ will be called {\em the reduced system} (\RS) in 3D.
\end{definition}

\begin{definition}\label{def:sol-RS-3D}\rm
We say that $\qnm$ satisfies \RS (or, informally, is a solution of \RS) iff
$$
\left\{
\begin{array}{rcl}
 \Fred_i(q_1,\dots,q_{n-1}) & = & 0, \quad i\in \{1,\dots,n-3\},\\[1ex]
 \hline
  \Fred_{n-2,x}(q_1,\dots,q_{n-1}) & = & 0, \\[1ex]
  \Fred_{n-2,y}(q_1,\dots,q_{n-1}) & = & 0, \\[1ex]
  z_{n-2} & = & 0,\\[1ex]
  \hline
  \Fred_{n-1,x}(q_1,\dots,q_{n-1}) & = & 0, \\[1ex]
  y_{n-1} & = & 0,\\[1ex]
  z_{n-1} & = & 0,\\[1ex]
  \hline
  q_n & = & q_n(q_1,\dots,q_{n-1}).
\end{array}
\right.
$$
\end{definition}

Note that \RS\  coincides with the system (\ref{eqn:reduced}), with the
equations for $y_{k_1},z_{k_1},z_{k_2}$ omitted.
Observe also that \RS\ no longer has $O(3)$ as a symmetry group. But still it is symmetric with respect to the reflections against
the coordinate planes.

The next theorem addresses the question: whether from \RS\ we
obtain the solution of (\ref{eq:cc-kart})?
\begin{theorem}\cite[Theorem 2]{MZ20}
\label{thm:red-to-full3D}
Assume that $\qnm$ is a solution of \RS\ satisfying
\begin{equation}
x_{n-1} \neq x_n. \label{eq:A1D3}
\end{equation}
\begin{description}
\item[Case 1]
   If  the vectors $(x_{n-1}-x_n, y_{n-1}-y_n)$ and $(x_{n-2}-x_n,y_{n-2}-y_n)$ are linearly independent,
     then
    $\qnm$ is a normalized central configuration, i.e.\ it satisfies (\ref{eq:cc-kart}).

\item[Case 2]    If  $\qnm$ is a solution such that $z_i=0$  for $i=1,\dots,n$, then $\qnm$  is a normalized central configuration, i.e. it satisfies (\ref{eq:cc-kart}).
\end{description}
\end{theorem}

Observe that condition appearing in case 1 of the above theorem is never satisfied for collinear solutions and also might not be satisfied for some planar solutions containing
three collinear bodies (such solutions exist  for $n=5$ and more, see \cite[Sec. A.2]{MZ}). This is why we included the second assertion in Theorem~\ref{thm:red-to-full3D}.

Another issue is how to determine whether a particular solution of the reduced system~(\ref{eq:cc-red-i}--\ref{eq:cc-red-yk2}) lies in the plane
$\{z=0\}$,  if we only know, that some multidimensional cube contains a unique solution of \RS. This issue is discussed in \cite{MZ20} in Section 3.3.

\begin{theorem}
\label{thm:nCC3D-RSsol}
Let $\overline{q} = (q_1, \ldots, q_{n-1})$ and $q = (\overline{q}, q_n(\overline{q}))$.
 Assume that  $\qnm = (q, m)$   is a  nCC. Then, in a suitable coordinate system and after  some permutation of bodies,  $\overline{q}$ is a solution of \RS\ satisfying $x_n \neq x_{n-1}$ and $x_{n-1} \neq 0$. Moreover, if
 $q$ is not collinear, then the vectors $(x_{n-1}-x_n, y_{n-1}-y_n)$ and $(x_{n-2}-x_n,y_{n-2}-y_n)$ are linearly independent
\end{theorem}
\proof
We can assume  that after a suitable permutation of bodies  $|q_{n-1}|$ is maximal, and in a suitable coordinate system $q_{n-1}=(x_{n-1},0,0)$ with $x_{n-1}>0$. From this it follows immediately that  $x_{n-1}>x_i$  for all $i \neq n-1$. From Theorem~\ref{thm:red-to-full3D} it follows that if $q$ is collinear, then it solves \RS.

In the non-collinear case we need to make further coordinate changes and permutations of bodies.  We look for the body is not on $OX$-axis. After a permutation (such that $(n-1) \mapsto (n-1)$) and suitable rotation we can assume that this is $(n-2)$-th body and $q_{n-2}=(x_{n-2},y_{n-2},0)$, $y_{n-2}>0$. Now we prove that there exists $j \neq n-1,n-2$ such that $(x_{n-1}-x_j, y_{n-1}-y_j)$ and $(x_{n-2}-x_j,y_{n-2}-y_j)$ are linearly independent.  Assume the contrary, then it for all $j$
$(x_j,y_j)$ must belong a line passing through $q_{n-1}$ and $q_{n-2}$ and the same is true for $(c_x,c_y)$ - a  projection on $OXY$ plane of the center of mass. The line connecting  $q_{n-2}$ and $q_{n-1}$
does not pass through the origin, but the center of mass is at the origin. So we obtain a contradiction. Therefore there exists $j$ such that  $(x_{n-1}-x_j, y_{n-1}-y_j)$ and $(x_{n-2}-x_j,y_{n-2}-y_j)$ are linearly independent. Now we permute bodies so that $j \mapsto n$, $n-1\mapsto n-1$ and $n-2\mapsto n-2$. From Theorem~\ref{thm:red-to-full3D} $q$ solves \RS.
\qed

\section{Non-degeneracy of CCs and the reduced systems }
\label{sec:non-degCC}
From the point of view of computer assisted proof (CAP)  the non-degeneracy plays a crucial role. If the solution of equations is non-degenerate then it is isolated and there is good chance to be verifiable  by a CAP.  In this section we will discuss the non-degeneracy of nCC's as solutions of \RS.   Some preliminary results in this direction in the planar case are contained in \cite{MZ,MZ22}. In the present paper  for $d=2$ and $d=3$ we will identify all situations, where the degeneracy might be result of passing from (\ref{eq:cc-kart}) to \RS. This is later  used by our program to avoid such situations - see Section~\ref{sec:about-prog}.

We begin with an adaption of  definition of non-degeneracy of CCs proposed by Moeckel  \cite[Def. 5]{Mlect2014}. 
The idea of Moeckel behind his notion of non-degeneracy is to allow only for degeneracy arising from the rotational symmetry of the problem.

Before we state our definition of non-degeneracy, let us notice that for any configuration $q$ the set $\mathcal{SO}(d)  q = \{ Rq \, | R \in \mathcal{SO}(d) \}$ is a smooth manifold, hence it makes sense to speak of dimension of $\mathcal{SO}(d)q$.
\begin{definition}\label{def:non-deg-cc}
  Assume $\qnm = (q, m)$ is a normalized central configuration, i.e.\ $F(q) = 0$, where $F$ is a system of equations~(\ref{eq:vector-field}). We say that $\qnm$ is \emph{non-degenerate}
  if
    $$rank(D\!F(q)) = dn-\dim \left(\mathcal{SO}(d)q\right).$$
  Otherwise the configuration is called \emph{degenerate}.
\end{definition}
If $d=2$, then $\dim \left(\mathcal{SO}(d)q\right)=1$ for any configuration $q$ without collision. For $d=3$ for configurations without collisions we have $\dim \left(\mathcal{SO}(d)q\right)=2$ for collinear configurations and $\dim \left(\mathcal{SO}(d)q\right)=3$ otherwise.

Now let us recall the standard definition of non-degeneracy of solution of a system of equations.
\begin{definition}\label{def:non-degeneracy}
Let $F\colon \reals^n\to\reals^n$ be a $\mathcal{C}^1$ function. A solution $x_0$ of $F(x) = 0$ is {\em non-degenerate} iff $DF(x_0)$ is an isomorphism.
\end{definition}

Note that under this definition no nCC's can be non-degenerate, because none is isolated. Note that there might be other reasons for the degeneracy of the solutions, for example being a bifurcation point of  nCC's as masses change.

\subsection{Center of mass reduction and the rank of Jacobian matrix }

\begin{lemma}\rm\label{lem:rank-com-red}
Let $\qnm$ be a  nCC. Then
\begin{equation*}
  \rank\left(D\! \Fred(q_1,\dots,q_{n-1})\right) = \rank\left(D\! F(q_1,\dots,q_{n-1},q_n)\right)-d. 
  \end{equation*}
\end{lemma}

\proof
This result follows from Lemma~\ref{lem:rank-after-elim} in Section~\ref{sssec:lm-rank}.
\qed

\subsubsection{Basic lemma about the Jacobian matrix}

 For $d = 3$ and $q = (q_1, \ldots, q_{n-1})$ let us denote:

\begin{eqnarray*}
D\Fred(q) & = &
\begin{bmatrix}
\smash[b]{\block{D\Fred_{1,x}(q)}{ \frac{\partial \Fred_{1,x}}{\partial x_1}(q)\ \ & \frac{\partial \Fred_{1,x}}{\partial y_1}(q)\ \ & \frac{\partial \Fred_{1,x}}{\partial z_1}(q)\  \ & \dots\ \ & \frac{\partial \Fred_{1,x}}{\partial y_{n-1}}(q)\ \ & \frac{\partial \Fred_{1,x}}{\partial z_{n-1}}(q)
}} \\[4em]
\smash[b]{\block{D\Fred_{1,y}(q)}{ \frac{\partial \Fred_{1,y}}{\partial x_1}(q)\ \  & \frac{\partial \Fred_{1,y}}{\partial y_1}(q)\ \  & \frac{\partial \Fred_{1,y}}{\partial z_1}(q)\ \  & \dots & \frac{\partial \Fred_{1,y}}{\partial y_{n-1}}(q)\ \  & \frac{\partial \Fred_{1,y}}{\partial z_{n-1}} (q)
}} \\[4em]
\smash[b]{\block{D\Fred_{1,z}(q)}{ \frac{\partial \Fred_{1,z}}{\partial x_1}(q)\ \  & \frac{\partial \Fred_{1,z}}{\partial y_1}(q)\ \  & \frac{\partial \Fred_{1,z}}{\partial z_1}(q)\ \  & \dots & \frac{\partial \Fred_{1,z}}{\partial y_{n-1}}(q)\ \  & \frac{\partial \Fred_{1,z}}{\partial z_{n-1}} (q)
}} \\[3em]
\ldots\\[1em]
\smash[b]{\block{D\Fred_{n-1,z}(q)}{ \frac{\partial \Fred_{n-1,z}}{\partial x_1}(q) & \frac{\partial \Fred_{n-1,z}}{\partial y_1}(q) & \frac{\partial \Fred_{n-1,z}}{\partial z_1}(q) & \dots & \frac{\partial \Fred_{n-1,z}}{\partial y_{n-1}}(q) & \frac{\partial \Fred_{n-1,z}}{\partial z_{n-1}} (q)
}}
\end{bmatrix}
\\[4em]
& = & \begin{bmatrix}
\mbox{}\quad  \frac{\partial \Fred}{\partial x_1}(q) \mbox{}\quad &   \frac{\partial \Fred}{\partial y_1}(q)
                                                         & \mbox{}\quad  \dots
                                                         & \mbox{}\quad \dots
                                                         &  \mbox{}\quad \frac{\partial \Fred}{\partial y_{n-1}}(q)
                                                         & \mbox{}\quad  \frac{\partial \Fred}{\partial z_{n-1}}(q)
                                                    \end{bmatrix} .
\end{eqnarray*}

Observe that  in the above matrix, rows are $D\Fred_{i,x}(q)$, $D\Fred_{i,y}(q)$ and  $D\Fred_{i,z}(q)$, while columns are $\frac{\partial \Fred}{\partial x_i}(q)$, $\frac{\partial \Fred}{\partial y_i}(q)$ and $\frac{\partial \Fred}{\partial z_i}(q)$.
In the case $d= 2$, the rows and columns corresponding to the
$z$-coordinate simply do not appear.

In the next lemma we show that some explicit linear combinations of rows or columns in $D\Fred(q)$ vanish if $q$ is nCC.
\begin{lemma}
\label{lem:DR-lin-komb-vanish3D}
Let  $(\overline{q},q_n({\overline{q}})) = \qnm$.
   Assume that  $\qnm$ is a  nCC.
   Then
   \begin{eqnarray}
      0&=&\sum_{i=1}^{n-1}m_i \left((x_i-x_n)  D \Fred_{i,y}(q)- (y_i-y_n)D \Fred_{i,x}(q) \right),  \label{eq:lkv-3d-xy} \\
   0&=&\sum_{i=1}^{n-1}m_i \left((y_i-y_n)  D \Fred_{i,z}(q)- (z_i-z_n)D \Fred_{i,y}(q) \right), \label{eq:lkv-3d-yz} \\
  0&=&\sum_{i=1}^{n-1}m_i \left((z_i-z_n)  D \Fred_{i,x}(q)- (x_i-x_n)D \Fred_{i,z}(q) \right), \label{eq:lkv-3d-zx} \\
   0&=&\sum_{i=1}^{n-1}\left( -\frac{\partial \Fred}{\partial x_i}(q)y_i + \frac{\partial \Fred}{\partial y_i}(q)x_i \right),  \label{eq:lkk-xy}\\
   0&=&\sum_{i=1}^{n-1}\left( -\frac{\partial \Fred}{\partial y_i}(q)z_i + \frac{\partial \Fred}{\partial z_i}(q)y_i \right),  \label{eq:lkk-yz} \\
    0&=&\sum_{i=1}^{n-1}\left( -\frac{\partial \Fred}{\partial z_i}(q)x_i + \frac{\partial \Fred}{\partial x_i}(q)z_i \right),\label{eq:lkk-zx}
   \end{eqnarray}
    where $q_i = (x_i,y_i,z_i)$ for all $i = 1, \ldots, n-1$.
\end{lemma}

\proof
From  Lemma~\ref{lem:c-angmom-red},  for arbitrary $(q_1,\dots,q_{n-1})$ with $q_n=q_n(q_1,\dots,q_{n-1})$ computed from the center of mass condition, we have
\begin{equation}
 0=\sum_{i=1}^{n-1}m_i (q_i-q_n) \wedge \Fred_i(q_1,\dots,q_{n-1}).
\end{equation}
By taking partial derivatives of the above equation with respect to
$x_j$, $y_j$ or $z_j$ (for $j=1,\dots,n-1$), and evaluating at $q$ (we have $\Fred_i(q)=0$), we obtain for $j=1,\dots,n-1$
\begin{eqnarray*}
  0&=&\sum_{i=1}^{n-1}m_i(q_i-q_n) \wedge \partial_v \Fred_i(q_1,\dots,q_{n-1}),
\end{eqnarray*}
 where $\partial_v = \frac{\partial }{\partial v}$ with $v \in \{x_j,y_j,z_j\}$.
Written component-wise, this gives us the following three equations:
\begin{eqnarray*}
  0&=& \sum_{i=1}^{n-1}m_i\left((x_i-x_n)  \partial_v \Fred_{i,y}(q)- (y_i-y_n)  \partial_v \Fred_{i,x}(q) \right), \\
  0&=& \sum_{i=1}^{n-1}m_i\left((y_i-y_n)  \partial_v \Fred_{i,z}(q)- (z_i-z_n)  \partial_v \Fred_{i,y}(q) \right), \\
  0&=& \sum_{i=1}^{n-1}m_i\left((z_i-z_n)  \partial_v \Fred_{i,x}(q)- (x_i-x_n)  \partial_v \Fred_{i,z}(q) \right).
\end{eqnarray*}
In terms of the rows of $D \Fred(q)$, the above equations can be written as (\ref{eq:lkv-3d-xy},\ref{eq:lkv-3d-yz},\ref{eq:lkv-3d-zx}), respectively.

Now we will establish (\ref{eq:lkk-xy},\ref{eq:lkk-yz},\ref{eq:lkk-zx}).
Let $O(t)$ be the rotation by angle $t$ in the $OXY$ plane.  It acts on configuration $q$ as follows: $q_i(t) = O(t)q_i$. Then $q(t)=(q_1(t),\dots,q_{n-1}(t),q_n(t))$
is nCC if $q$ is, that is,
 \begin{equation}
  \Fred(O(t)q)=0. \label{eq:R0q3d}
\end{equation}

   Observe that $\frac{d}{dt}q_i(t)_{t=0}=(-y_i,x_i)$.
By  taking the derivative of (\ref{eq:R0q3d}) with respect to the angle  for $t=0$ we obtain
\begin{equation*}
  0=\sum_{i=1}^{n-1}\left( -\frac{\partial \Fred}{\partial x_i}(q)y_i + \frac{\partial \Fred}{\partial y_i}(q)x_i \right).
\end{equation*}
This is equation (\ref{eq:lkk-xy}). Equations (\ref{eq:lkk-yz},\ref{eq:lkk-zx}) are obtained by
rotations in other coordinate planes.
\qed

Now we are ready to establish relation between non-degenerate nCCs (in the sense of Def.~\ref{def:non-deg-cc}) and non-degenerate solutions of \RS.
In the sequel, by $D\!\mbox{\RS}(q)$ we will denote the Jacobian matrix of \RS~at~$q$.

\subsection{Non-degeneracy of  normalized CCs  in 2D }\label{subsec:ndeg2d}
The goal of this section is to discuss the relation between non-degenerate nCC's and non-degenerate solutions of \RS\ in  the planar case. We give a complete description of situations in which the degeneracy is produced when passing to \RS. This insight allows us to avoid this phenomenon when realizing a computer assisted proof.

In the sequel, first, we state the theorems in the Results section, and subsequently provide their proofs in the Proofs section.

\subsubsection{Results}
Let $\overline{q} = (q_1, \ldots, q_{n-1})$, $q = (\overline{q}, q_n(\overline{q}))$ and $\qnm = (q, m)$.

\begin{theorem}
\label{thm:RS-nondeg-equiv}
 Assume that  $\qnm$  is an  nCC with  $y_{n-1}=0$ such that
 \begin{equation}
  x_{n-1} \neq x_n, \quad x_{n-1}\neq 0.  \label{eq:2D-nondegproj}
 \end{equation}
 Then $\qnm$ is a non-degenerate nCC iff $\qnm$ is a non-degenerate solution of \RS.
\end{theorem}
The implication  $\Leftarrow$ has been proved  in \cite[Theorem 14]{MZ} and in \cite[Theorem 7]{MZ22} .

\begin{theorem}\cite[Theorem 5]{MZ22}
\label{thm:sol-RS-deg}
 Assume that  $\qnm $   is an  nCC such that $y_{n-1}=0$.
  If $x_{n-1} = x_n$ or $x_{n-1}=0$, then $\qnm$ is a degenerate solution
of \RS.
\end{theorem}

Theorems~\ref{thm:RS-nondeg-equiv} and~\ref{thm:sol-RS-deg} together with Theorem~\ref{thm:non-deg-RS} show that after a suitable rotation and permutation of bodies  each non-degenerate nCC is a non-degenerate solution
of \RS.  Of course the suitable rotation and permutation depends on nCC, but in fact we just need to place one selected body, which is not at the origin, on $OX$-axis and
do any permutation after which it becomes $(n-1)$-th body.  This is realized in our program, see Section~\ref{sec:about-prog}.

\subsubsection{Proofs}
To prove Theorem~\ref{thm:RS-nondeg-equiv}, first let us state the following lemma:

 \begin{lemma}\cite[Lemma 4]{MZ22}
\label{lem:rankRS}
 Assume that   $(q, m)$  is an  nCC with $y_{n-1}=0$  satisfying
 \begin{equation*}
 x_n \neq x_{n-1},  \quad x_{n-1} \neq 0.
 \end{equation*}
 Then $\rank(D\!\mbox{\RS}(q)) = \rank(D\! \Fred(q))$.
\end{lemma}

\proof (of Lemma~\ref{lem:rankRS})
Observe that the jacobian matrix of \RS\ is equal to  matrix $D\! \Fred(q)$ from which we removed the last column (which is the consequence of the restriction to $y_{n-1}=0$) and the last row (which is the consequence of dropping the equation $\Fred_{n-1,y}(q)=0$). We need to show that such removal does not change the rank of the matrix.

Equation (\ref{eq:lkv-3d-xy}) in Lemma~\ref{lem:DR-lin-komb-vanish3D} shows that, if $x_{n-1}-x_n \neq 0$, then   the last row (i.e. $D\! \Fred_{n-1,y}$)  can be expressed as the linear combination of other rows, hence it can be removed from the matrix without changing its rank.

Equation (\ref{eq:lkk-xy}) in Lemma~\ref{lem:DR-lin-komb-vanish3D} shows, that
if $x_{n-1} \neq 0$, then we can express  $\frac{\partial \Fred}{\partial y_{n-1}}$ (the last column in the matrix $D \Fred$)
in terms of other columns.

Therefore we can remove the last row and the last column from the matrix $D \Fred(q)$ without decreasing its rank.
\qed

\proof (of Theorem~\ref{thm:RS-nondeg-equiv}) From Lemmas~\ref{lem:rank-com-red} and Lemma~\ref{lem:rankRS} it follows that
\begin{equation}
  \rank(D\! \mbox{\RS}(q) )=\rank\left(D\! \Fred(q)\right)=\rank \left(D\! F(q)\right)-2.
\end{equation}

Therefore $\rank( D\! \mbox{\RS}(q) )=2n-3$ iff $\rank \left(D\! F(q)\right)=2n-1$.  This finishes the proof.
\qed

\noindent
\proof  (of Theorem~\ref{thm:sol-RS-deg})
We use Lemma~\ref{lem:rankRS}.
If $x_{n-1} = x_n$, then equation (\ref{eq:lkv-3d-xy}) in Lemma~\ref{lem:DR-lin-komb-vanish3D} gives a vanishing linear combination of rows in $D\!\mbox{\RS}(q)$, because $x_{n-1}-x_n$ multiplying row $D\Fred_{n-1,y}$ (which is not a row in $D\!\mbox{\RS}(q)$) vanishes.
Observe that some coefficients  in this linear combinations must be non-zero, otherwise we will have $q_i=q_n$ for all $i$.

If $x_{n-1}=0$, then  equation (\ref{eq:lkk-xy})  in Lemma~\ref{lem:DR-lin-komb-vanish3D}  gives a vanishing linear combination of columns of the jacobian matrix of \RS\ at $q$, because the coefficient multiplying column $\frac{\partial \Fred}{\partial y_{n-1}}$ (which is not a column in $D\!\mbox{\RS}(q)$) vanishes.
Observe that some coefficients  in this linear combinations must be non-zero, otherwise we will have $q_i=0$ for all $i$.

Hence in both cases the rank of the jacobian matrix of \RS\ at $q$ cannot be maximal.
\qed

\subsection{Non-degeneracy of   normalized CCs  in 3D}
\label{subsec:ndeg3d}

The goal of this section is to  give a complete description of situations in which the degeneracy of solutions of \RS\ for non-degenerate nCC's  is produced when passing to \RS\ for $d=3$.
Compared to the 2D case, in 3D there are more possibilities that can lead to degeneracies in the solutions of \RS.
 For instance, if $y_{n-2} = 0$, then rotating a solution of \RS\ around the OX axis yields a continuous family (a circle) of solutions of \RS, unless the configuration is collinear.

\subsubsection{Results}

\begin{theorem}\label{thm:RS3D-nondeg-equiv}
 Assume that   $\qnm$  is an nCC satisfying
 the normalization  $y_{n-1}=z_{n-1}=z_{n-2}=0$.
 If
 \begin{eqnarray}
   x_{n-1}-x_n \neq 0,\  x_{n-1} \neq 0,\ y_{n-2} & \neq & 0   \label{eq:3Dnondeg-c1} \\
  \det \left[\begin{array}{cc}
  (y_{n-2}-y_n), &   (y_{n-1}-y_n)  \\
  (x_{n-2}-x_n),  & (x_{n-1}-x_n)  \\
\end{array}
\right] &\neq & 0,  \label{eq:3Dnondeg-c2}
 \end{eqnarray}
 then $\qnm$ is non-degenerate nCC iff $\qnm$ is a non-degenerate solution of \RS.
\end{theorem}

\begin{theorem}
\label{thm:sol-RS3D-deg}
 Assume that   $\qnm$  is  an nCC  satisfying
 the normalization $y_{n-1}=z_{n-1}=z_{n-2}=0$.
\begin{description}
\item[Case 1] If $x_{n-1} = x_n$ or $x_{n-1}=0$,  then $\qnm$ is a degenerate solution
of \RS.

\item[Case 2] If $y_{n-2}=0$ and $\qnm$ is  non-collinear,
 then $\qnm$ is a degenerate solution
of \RS.

\item[Case 3] If
\begin{equation}
 \det \left[\begin{array}{cc}
  (y_{n-2}-y_n) &   (y_{n-1}-y_n)  \\
  (x_{n-2}-x_n)  & (x_{n-1}-x_n)  \\
\end{array}
\right] = 0  \label{eq:det=0-1}
\end{equation}
and $\qnm$ is non-collinear, then $\qnm$ is degenerate solution of \RS.
\end{description}
\end{theorem}

 The above theorems provide a complete characterization of the situations in which degeneracy arises when passing to the reduced system \RS, assuming that the $n$–body central configuration under consideration is non-collinear. Theorem~\ref{thm:nCC3D-RSsol}—in fact, its proof—explains how to choose appropriate rotations of the coordinate system and permutations of the bodies in order to avoid such degeneracies in the non-collinear case.

Collinear central configurations require separate treatment. Although collinear $n$–body central configurations are non-degenerate as central configurations, they may nevertheless appear as degenerate solutions of the reduced system \RS. This issue is discussed in Section~\ref{sec:collncc}.


\subsubsection{Proofs}
 \begin{lemma}\label{lem:rankRS3D}
 Assume that  nCC $\qnm$     satisfies
 the normalization  $y_{n-1}= z_{n-1} = z_{n-2}=0$.

 Assume that
 \begin{eqnarray}
 x_{n-1}-x_n \neq 0, \quad x_{n-1} \neq 0, \quad  y_{n-2} & \neq & 0, \\
  \det \left[\begin{array}{cc}
  (y_{n-2}-y_n), &   (y_{n-1}-y_n)  \\
  (x_{n-2}-x_n),  & (x_{n-1}-x_n)  \\
\end{array}
\right] & \neq & 0.  \label{eq:det-neq-0}
\end{eqnarray}

 Then the rank of $D\! \mbox{\RS}(q)$ is equal  to the rank of $D\Fred(q)$.
\end{lemma}

\proof (of Lemma~\ref{lem:rankRS3D})
 Observe that the Jacobian matrix of \RS\ is obtained
from matrix $D\! \Fred$ by removing  three  columns, corresponding to $\frac{\partial \Fred}{\partial y_{n-1}}$, $\frac{\partial \Fred}{\partial z_{n-1}}$ and $\frac{\partial \Fred}{\partial z_{n-2}}$ (due to restriction to $y_{n-1}=z_{n-1}=z_{n-2}=0$), and  three  rows, corresponding to $D \Fred_{n-1,y}$, $D \Fred_{n-1,z}$ and $D \Fred_{n-2,z}$ (due to dropping  equations $\Fred_{n-1,y}=0$, $\Fred_{n-1,z}=0$,
$\Fred_{n-2,z}=0$). We need to show that this removal does not change the rank of the matrix.

First we want to remove  rows $D \Fred_{n-1,y}$, $D \Fred_{n-1,z}$, $D \Fred_{n-2,z}$.
\begin{itemize}
\item From (\ref{eq:lkv-3d-xy}) we obtain
\begin{eqnarray*}
  -m_{n-1}(x_{n-1}-x_n)  D \Fred_{n-1,y}=\sum_{i=1}^{n-2}m_i(x_i-x_n)  D \Fred_{i,y}- \sum_{i=1}^{n-1}m_i(y_i-y_n)D \Fred_{i,x}.
\end{eqnarray*}
We see that if $x_{n} \neq x_{n-1}$, then $ D \Fred_{n-1,y}$ can be expressed as a linear combination of certain rows from $D\!\mbox{\RS}$.

\item From (\ref{eq:lkv-3d-yz}) we obtain
\begin{eqnarray*}
 -m_{n-1}(y_{n-1}-y_n)  D \Fred_{n-1,z} &-& m_{n-2}(y_{n-2}-y_n)  D \Fred_{n-2,z}  =  \sum_{i=1}^{n-3}m_i(y_i-y_n)  D \Fred_{i,z} \\
& & - \sum_{i=1}^{n-2}m_i(z_i-z_n)D \Fred_{i,y}  - (z_{n-1}-z_n)m_{n-1}D \Fred_{n-1,y}
\end{eqnarray*}
Observe that rows on the rhs can be expressed by rows from  $D\!\mbox{\RS}$.

\item From (\ref{eq:lkv-3d-zx}) we obtain
\begin{eqnarray*}
   m_{n-1}(x_{n-1}-x_n)D \Fred_{n-1,z}   + m_{n-2}(x_{n-2}-x_n)D \Fred_{n-2,z} & = & \sum_{i=1}^{n-1}m_i(z_i-z_n)  D \Fred_{i,x}\\
    & & - \sum_{i=1}^{n-3}m_i(x_i-x_n)D \Fred_{i,z}.
\end{eqnarray*}
\end{itemize}
We obtain a system of two equations in which the left-hand side consists of a linear combination of the rows
$D \Fred_{n-1,z}$ and $D \Fred_{n-2,z}$, and the right-hand side involves other rows from  $D\!\mbox{\RS}$. Therefore, it suffices for the determinant of the coefficient matrix on the left-hand side to be nonzero in order to express the rows to be removed as linear combinations of rows from  $D\!\mbox{\RS}$. This determinant is nonzero by assumption~(\ref{eq:det-neq-0}).

Now we want to remove  columns $\frac{\partial \Fred}{\partial y_{n-1}}$, $\frac{\partial \Fred}{\partial z_{n-1}}$ and $\frac{\partial \Fred}{\partial z_{n-2}}$.
\begin{itemize}
\item From (\ref{eq:lkk-xy}) we obtain
\begin{eqnarray*}
   -x_{n-1}\frac{\partial \Fred}{\partial y_{n-1}} & = & -\sum_{i=1}^{n-1}\frac{\partial \Fred}{\partial x_i}y_i + \sum_{i=1}^{n-2}\frac{\partial \Fred}{\partial y_i}x_i.
\end{eqnarray*}
We see that if $ x_{n-1} \neq 0$, then we can express $\frac{\partial \Fred}{\partial y_{n-1}}$ as the linear combination of some columns appearing in  $D\!\mbox{\RS}$ .

\item From (\ref{eq:lkk-yz}) we obtain
 \begin{eqnarray*}
 -\frac{\partial \Fred}{\partial z_{n-1}}y_{n-1}  - \frac{\partial \Fred}{\partial z_{n-2}}y_{n-2}  & = & -\sum_{i=1}^{n-1}\frac{\partial \Fred}{\partial y_i}z_i + \sum_{i=1}^{n-3}\frac{\partial \Fred}{\partial z_i}y_i.
 \end{eqnarray*}
Due to normalization, we have $y_{n-1}=0$ and $z_{n-1}=0$, and the above equation  becomes
 \begin{eqnarray*}
   - \frac{\partial \Fred}{\partial z_{n-2}}y_{n-2}  & = & -\sum_{i=1}^{n-2}\frac{\partial \Fred}{\partial y_i}z_i + \sum_{i=1}^{n-3}\frac{\partial \Fred}{\partial z_i}y_i.
 \end{eqnarray*}
 By the assumption $y_{n-2} \neq 0$, hence column $\frac{\partial \Fred}{\partial z_{n-2}}$ can be expressed by columns from  $D\!\mbox{\RS}$.

\item  From (\ref{eq:lkk-zx}) we obtain
 \begin{eqnarray}
    \frac{\partial \Fred}{\partial z_{n-1}}x_{n-1} + \frac{\partial \Fred}{\partial z_{n-2}}x_{n-2} & = & -\sum_{i=1}^{n-3}\frac{\partial \Fred}{\partial z_i}x_i + \sum_{i=1}^{n-1}\frac{\partial \Fred}{\partial x_i}z_i.
 \end{eqnarray}
 \end{itemize}
 Given the assumption $x_{n-1} \neq 0$, the column $\frac{\partial \Fred}{\partial z_{n-1}}$ can be written as a linear combination of columns in  $D\!\mbox{\RS}$ and $\frac{\partial \Fred}{\partial z_{n-2}}$, which, as previously shown, is itself expressible in terms of columns of  $D\!\mbox{\RS}$.
\qed

\proof (of Theorem~\ref{thm:RS3D-nondeg-equiv})
Observe that condition (\ref{eq:3Dnondeg-c2}) implies that $\qnm$ is not collinear.  Therefore it is enough to prove that
\begin{equation}
 \rank(D\! \mbox{\RS\ }(q) )=3n-6 \quad \mbox{iff} \quad \rank \left(D\! F(q)\right)=3n-3. \label{eq:ranks-deg3dpr}
\end{equation}

From Lemmas~\ref{lem:rank-com-red} and Lemma~\ref{lem:rankRS3D} it follows that
\begin{equation}
  \rank(D\! \mbox{\RS\ }(q) )=\rank\left(D\! \Fred(q)\right)=\rank \left(D\! F(q)\right)-3.
\end{equation}
This establishes (\ref{eq:ranks-deg3dpr}) and finishes the proof.
\qed

\noindent
\proof (of Theorem~\ref{thm:sol-RS3D-deg})
\begin{description}
\item[Case 1]  The proof  is the same as  in 2D case; see Theorem~\ref{thm:sol-RS-deg} and its proof.

\item[Case 2] If $y_{n-2}=0$, then  we can rotate a solution of \RS\
 around OX axis to obtain a whole circle of solutions of \RS\ (unless the configuration is collinear). Hence, the configuration as a solution of \RS\ is degenerate.

\item[Case 3] Our goal is to find a non-trivial vanishing linear combination of rows in matrix $D\! \mbox{\RS}(q)$.
From previous reasoning we can assume that $x_{n-1}-x_n \neq 0$, because otherwise $\qnm$ is a degenerate solution of \RS.
Our point of departure are the  identities  for rows in $D \Fred$ from Lemma~\ref{lem:DR-lin-komb-vanish3D}
\begin{eqnarray}
  -(x_{n-1}-x_n) m_{n-1} D \Fred_{n-1,y} & = & \sum_{i=1}^{n-2}(x_i-x_n) m_i D \Fred_{i,y}\nonumber\\
  & & - \sum_{i=1}^{n-1}(y_i-y_n)m_i D \Fred_{i,x}, \label{eq:DRn1y}\\
 -(y_{n-1}-y_n) m_{n-1} D \Fred_{n-1,z} -(y_{n-2}-y_n) m_{n-2} D \Fred_{n-2,z} & = & \sum_{i=1}^{n-3}(y_i-y_n) m_i D \Fred_{i,z} \nonumber  \\
  - \sum_{i=1}^{n-2}(z_i-z_n)m_i D \Fred_{i,y}
 - (z_{n-1}-z_n)m_{n-1}D \Fred_{n-1,y}, \label{eq:DRn12nd} \\
   (x_{n-1}-x_n)m_{n-1}D \Fred_{n-1,z}   + (x_{n-2}-x_n)m_{n-2}D \Fred_{n-2,z} & = & \sum_{i=1}^{n-1}(z_i-z_n) m_i D \Fred_{i,x} \nonumber \\
 & &    - \sum_{i=1}^{n-3}(x_i-x_n)m_iD \Fred_{i,z}. \label{eq:DRn13nd}
\end{eqnarray}
Observe that on lhs of above identities we have only rows which are not present in  $D\! \mbox{\RS}(q)$ and on rhs rows from $D\! \mbox{\RS}(q)$
and single row $D \Fred_{n-1,y}$, which is not from $D\! \mbox{\RS}(q)$, but which using (\ref{eq:DRn1y}) could be expressed as linear combination of rows in $D\! \mbox{\RS}(q)$.

From (\ref{eq:det=0-1}) it follows  that there exist $a,b \in \mathbb{R}$ with at least one of them non-zero such that
\begin{equation}
  b (x_{n-2}-x_n, x_{n-1}-x_n) - a (y_{n-2}-y_n,   y_{n-1}-y_n)=0.  \label{eq:abxyn1n2}
\end{equation}

Then  $a$ times  equation (\ref{eq:DRn12nd}) plus $b$ times of equation (\ref{eq:DRn13nd})
gives  equation with vanishing lhs. We obtain
\begin{eqnarray*}
 0 & = & \sum_{i=1}^{n-3}a(y_i-y_n)  m_i D \Fred_{i,z}  - \sum_{i=1}^{n-2}a(z_i-z_n)m_iD \Fred_{i,y}  - a(z_{n-1}-z_n)m_{n-1}D \Fred_{n-1,y} \\
& &  + \sum_{i=1}^{n-1}b(z_i-z_n) m_i D \Fred_{i,x}- \sum_{i=1}^{n-3}b(x_i-x_n)m_iD \Fred_{i,z},
\end{eqnarray*}
and further
\begin{eqnarray*}
 0 & = & \sum_{i=1}^{n-3}\left(a(y_i-y_n) - b(x_i-x_n)  \right)m_iD \Fred_{i,z}  - \sum_{i=1}^{n-2}a(z_i-z_n)m_iD \Fred_{i,y} \\
 & & - a(z_{n-1}-z_n)m_{n-1}D \Fred_{n-1,y}
  + \sum_{i=1}^{n-1}b(z_i-z_n) m_i D \Fred_{i,x}
\end{eqnarray*}
From (\ref{eq:DRn1y}) we can compute $m_{n-1}D \Fred_{n-1,y}$ and insert it in the above equation to obtain
\begin{eqnarray*}
 0 & = & \sum_{i=1}^{n-3}\left(a(y_i-y_n) - b(x_i-x_n)  \right)m_i D \Fred_{i,z}  - \sum_{i=1}^{n-2}a(z_i-z_n)m_iD \Fred_{i,y} \\
& &  + a\frac{z_{n-1}-z_n}{x_{n-1}-x_n}\left(\sum_{i=1}^{n-2}(x_i-x_n)  m_iD \Fred_{i,y}- \sum_{i=1}^{n-1}(y_i-y_n)m_iD \Fred_{i,x} \right) \\
& &  + \sum_{i=1}^{n-1}b(z_i-z_n) m_i D \Fred_{i,x}.
\end{eqnarray*}
Finally after regrouping we obtain
\begin{eqnarray}\label{eq:3d-vani-com-lin}
 0 & = & \sum_{i=1}^{n-3}\left(a(y_i-y_n) - b(x_i-x_n)  \right)m_iD \Fred_{i,z} \\
& &   + a\sum_{i=1}^{n-2}\left(\frac{z_{n-1}-z_n}{x_{n-1}-x_n}(x_i-x_n) -  (z_i-z_n)\right)m_iD \Fred_{i,y}  \nonumber \\
& &  + \sum_{i=1}^{n-1}\left(b(z_i-z_n) -  a\frac{z_{n-1}-z_n}{x_{n-1}-x_n}(y_i-y_n) \right) m_i D \Fred_{i,x}. \nonumber
\end{eqnarray}
We obtained a vanishing linear combination of rows in $D\! \mbox{\RS}(q)$.   We will prove that if all coefficients of~(\ref{eq:3d-vani-com-lin}) vanish, then
$q_i$'s are collinear. This will finish the proof.

Vanishing coefficients in the first sum will give   $-a(y_i-y_n) + b(x_i-x_n)=0$ for $i=1,\dots,n-3$. This, together with~(\ref{eq:abxyn1n2}) implies that
\begin{equation}
  b(x_1-x_n,x_2-x_n,\dots,x_{n-1}-x_n)-a(y_1-y_n, y_2 - y_n,\dots,y_{n-1}-y_n) =  0.  \label{eq:sum1=0}
\end{equation}
Vanishing of coefficients in second sum implies that
\begin{equation}
  \frac{z_{n-1}-z_n}{x_{n-1}-x_n}(x_1-x_n,\dots,x_{n-2}-x_n,x_{n-1}-x_n) - (z_1-z_n,\dots, z_{n-2}-z_n,z_{n-1}-z_n)=0. \label{eq:sum2=0}
\end{equation}
In the above equality the vanishing on the last coordinate is a consequence of the identity $\frac{z_{n-1}-z_n}{x_{n-1}-x_n}(x_{n-1}-x_n) - (z_{n-1}-z_n)$.

Finally, vanishing of coefficients in third sum implies that
\begin{equation}
b(z_1-z_n,\dots,z_{n-1}-z_n) - a\frac{z_{n-1}-z_n}{x_{n-1}-x_n}(y_1-y_n, y_2 - y_n,\dots,y_{n-1}-y_n)=0  \label{eq:sum3=0}
\end{equation}
Observe that identity (\ref{eq:sum3=0}) is not really an independent condition, because it follows from (\ref{eq:sum1=0},\ref{eq:sum2=0}).

Consider a matrix $M \in \mathbb{R}^{3 \times (n-1)}$ with $n-1$ columns given by $q_i-q_n$ for $i=1,\dots,n-1$, i.e.
$$  M=\begin{bmatrix}
    x_1-x_n, & x_2-x_n, & \dots & x_{n-1}-x_n \\
    y_1-y_n, & y_2-y_n, & \dots & y_{n-1}-y_n \\
     z_1-z_n, & z_2-z_n, & \dots & z_{n-1}-z_n
  \end{bmatrix}
$$
 From the above considerations, it follows that
if all the coefficients in~(\ref{eq:3d-vani-com-lin}) vanish, then all three rows of the matrix $M$ are proportional. Therefore, $M$ has rank equal to $1$.  Consequently, all columns are proportional,
that is, they are collinear, and we have
\begin{equation}
  q_i - q_{n}= \lambda_i (v_1,v_2,v_3)
\end{equation}
for some nonzero constants $\lambda_i$ and a nonzero vector  $(v_1,v_2,v_3)$.
We want to show that $q_i$'s  are collinear.
We have
\begin{equation}
  q_i=q_n + \lambda_i v  \label{eq:qiqnv}
\end{equation}
and we insert (\ref{eq:qiqnv}) into center of mass condition (\ref{eq:cc-cofmass-l}) to obtain
\begin{eqnarray*}
  0 & = & \sum_{i=1}^{n-1}m_iq_i +  m_{n}q_{n} = \sum_{i=1}^{n-1}m_i(q_n + \lambda_i v) + m_{n}q_{n}= \left(\sum_{i=1}^n m_i\right)q_n +
  \left(\sum_{i=1}^{n-1}m_i\lambda_i\right) v,
\end{eqnarray*}
hence
\begin{equation}
 q_n = - \left(\sum_{i=1}^n m_i\right)^{-1}\left(\sum_{i=1}^{n-1}m_i\lambda_i\right) v.
\end{equation}
Therefore the configuration is collinear.
\end{description}
\qed

\section{Collinear nCCs and the degeneracy issues}
\label{sec:collncc}

In the spatial case, in the context of the non-degeneracy of nCCs as solutions of  \RS\ 
 (see Theorems~\ref {thm:RS3D-nondeg-equiv} and~\ref {thm:sol-RS3D-deg}), collinear nCCs
 play a special role. Moreover, since it is known that such configurations exist (see Theorem~\ref{thm:mou}), it is natural 
 to ask whether they are non-degenerate solutions of \RS.  In the present section, we show that this is indeed the case for $d = 2$ (see Theorem 20), whereas the question for $d = 3$ remains open. We first recall two fundamental results concerning the existence of collinear central configurations. The first is Moulton’s theorem on existence.
 
 \begin{theorem}\cite{Mou}\cite[Proposition 18]{Mlect2014}
\label{thm:mou}
For any ordering of masses there exists a unique (up to scaling and translations) collinear central configuration.
Hence there are $n!/2$ classes of collinear central configurations.
\end{theorem}

The next one is a    theorem of Conley (see \cite[Prop. 19]{Mlect2014} or \cite{P}) 
\begin{theorem}
\label{thm:collcc-nondeg}
For any $d$, all normalized collinear central configurations  satisfy $\rank D\! F(q) = dn -(d-1)$.
\end{theorem}
As a consequence of the above result we obtain that all collinear nCCs are non-degenerate in the sense of Definition~\ref{def:non-deg-cc}.
In this section we are interested, whether for $d=2$ and $d=3$ a collinear nCC with $x_{n-1} \neq 0$ is non-degenerate as a solution of \RS. We will show that for $d=2$ this is true, but  for $d=3$ it may fail.

To see what happens when we pass to \RS\ , we need to have a  detailed knowledge of $D\! F(q)$. Therefore we need first go over the proof
of Theorem~\ref{thm:collcc-nondeg} formulating some key steps from its proof as lemmas, which will be later used to analise $D\!\mbox{\RS}(q)$.

We will use $R$ rather than $F$, since most of the arguments are formulated for $R$, owing to the fact that $DR$ is symmetric. This means that our equations for nCC are
\begin{equation}
  R_i(q_1,\dots,q_n)=0, \quad i=1,\dots,n.
\end{equation}
Throughout the remainder of this section we assume that
$q$ is a collinear, normalized CC lying on the $OX$-axis.


We are interested in the structure of $D\! R(q)$.  Let  $A \in \mathbb{R}^{n \times n}$ be given by
\begin{eqnarray*}
  A_{ii}&=&  \sum_{j,j\neq i}  \frac{m_i m_j}{r_{ij}^3}, \\
  A_{ij}&=& -\frac{m_i m_j}{r_{ij}^3}, \quad i \neq j.
\end{eqnarray*}
Observe that $A$ is symmetric. Let $M \in  \mathbb{R}^{n \times n}$ be a diagonal matrix with $M_{ii}=m_i$.
If we order variables as follows $(x_1,x_2,\dots,x_n,y_1,\dots,y_n,z_1,\dots,z_n)$, then it is easy to see that $D\! R(q)$ has a block diagonal structure: for each variable $x,y,z$ we have a block on diagonal,
\begin{equation}
  DR = \left[ \begin{array}{ccc}
                 M+2A & 0 & 0 \\
                 0 & M-A & 0 \\
                 0 & 0  & M-A
              \end{array}
      \right]
\end{equation}

The application of the Gershogorin  theorem gives the following result
\begin{lemma}
\label{lem:x-block-nondeg}
  The block in $x$-direction, i.e. $M+2A$, is positive definite.
\end{lemma}

The situation with block for $y$-variable (and $z$-variable) is more subtle.
The following statement can be found in  \cite{Mlect2014,P}  
\begin{lemma}
\label{lem:yblock-ev}
The matrix $M - A$ has exactly one positive eigenvalue, one zero eigenvalue, and $n-2$ negative eigenvalues.
\end{lemma}

Now we consider the center of mass reduction.
Let  $A^{\mathrm{red}} \in \mathbb{R}^{(n-1) \times (n-1)}$ be given by
\begin{eqnarray*}
  A^{\mathrm{red}}_{ii}&=&  \sum_{j=1, j \neq i}^{n-1}\frac{m_i m_j}{r_{ij}^3} + \frac{m_i (m_i + m_n)}{r_{in}^3} , \\
  A^{\mathrm{red}}_{ij}&=&  m_i m_j\left( \frac{1}{r_{in}^3}-\frac{1}{r_{ij}^3}\right), \quad i \neq j.
\end{eqnarray*}
Observe that $A^{\mathrm{red}}$ is not symmetric. Let $M^{\mathrm{red}}\in \mathbb{R}^{(n-1) \times (n-1)}$ be a diagonal matrix with $M^{\mathrm{red}}_{ii}=m_i$ with $i=1,\dots,n-1$.  Then we have
\begin{equation}
  DR^{\mathrm{red}} = \left[ \begin{array}{ccc}
                 M^{\mathrm{red}}+2A^{\mathrm{red}} & 0 & 0 \\
                 0 & M^{\mathrm{red}}-A^{\mathrm{red}} & 0 \\
                 0 & 0  & M^{\mathrm{red}}-A^{\mathrm{red}}
              \end{array}
      \right]
\end{equation}

We introduce the following notation for blocks
\begin{eqnarray*}
   DR^{\mathrm{red}}_x = \left[\frac{\partial R^{\mathrm{red}}_{i,x}}{\partial x_j}\right]_{i,j=1,\dots,n-1}, \\
    DR^{\mathrm{red}}_y = \left[\frac{\partial R^{\mathrm{red}}_{i,y}}{\partial y_j}\right]_{i,j=1,\dots,n-1}, \\
      DR^{\mathrm{red}}_z = \left[\frac{\partial R^{\mathrm{red}}_{i,z}}{\partial z_j}\right]_{i,j=1,\dots,n-1},
\end{eqnarray*}
and analogously (taking into account the variables and equations we drop) we define $D\mbox{\RS}_x$, $D\mbox{\RS}_y$ and $D\mbox{\RS}_z$.
Observe that for collinear central configurations holds
\begin{equation}
DR^{\mathrm{red}}_y=DR^{\mathrm{red}}_z= M^{\mathrm{red}}-A^{\mathrm{red}}.  \label{eq:DRredxy}
\end{equation}

For blocks $D\mbox{\RS}_v$,   $DR^{\mathrm{red}}_v$ for $v \in \{x,y,z\}$ we will apply the standard notion of (non-)degeneracy, i.e. the linear map represented by a given block should be
isomorphism to be  non-degenerate. 

\begin{lemma}
\label{lem:xblk-non-deg}
Assume that $q$ is normalized collinear central configuration contained in $OX$-axis. Then $D\!\Rred_x(q)=D\! \mbox{\RS}_x(q)$ is non-degenerate.
\end{lemma}
\proof
This follows from Lemma~\ref{lem:rank-com-red} for $d=1$ and Lemma~\ref{lem:x-block-nondeg}.
\qed

\begin{lemma}
\label{lem:cblk-y-nondeg}
Assume  that $x_{n-1} \neq 0$. Then
 $D\!\mbox{\RS}_y$ is non-degenerate.
\end{lemma}
\proof
For collinear nCCs $x_i \neq x_j$ for all $i \neq j$. First observe that
$D\!\mbox{\RS}_y$  does not depend on whether we work in dimension $d=2$ or $d=3$. Thus, we may restrict our attention to the case $d=2$, where Lemma~\ref{lem:rankRS} applies:
\begin{eqnarray*}
   \rank(D\! \mbox{\RS}(q))=\rank(D\! \Fred(q))=\rank(D\! \Rred(q)).
\end{eqnarray*}
From the block diagonal structure of all matrices involved and from Lemma~\ref{lem:xblk-non-deg} we obtain
\begin{equation*}
   \rank(D\! \mbox{\RS}_y(q))=\rank (DR^{\mathrm{red}}_y(q)).
\end{equation*}
We will compute $\rank (DR^{\mathrm{red}}_y(q))$.
 From Lemma~\ref{lem:rank-com-red} and the block-diagonal structure of $D\! R(q)$ and $D\! \Rred (q)$ we have
\begin{eqnarray*}
  \rank (D\! \Rred_x(q)) + \rank (D\!\Rred_y(q)) =  \rank (DR^{\mathrm{red}}(q))=\rank(DR(q))-2\\
  =\rank (D\! R_x(q)) + \rank (D\! R_y(q))-2
\end{eqnarray*}
hence
\begin{equation*}
 \rank (D\!\Rred_y(q))= \rank (D\! R_x(q)) - \rank (D\! \Rred_x(q)) + \rank (D\! R_y(q))-2.
\end{equation*}
From Lemma~\ref{lem:rank-com-red} it follows that $\rank (D\! R_x(q)) - \rank (D\! \Rred_x(q))=1$ and from Lemma~\ref{lem:yblock-ev}
$ \rank(DR_y(q))=n-1$, hence we obtain
\begin{equation*}
  \rank (DR^{\mathrm{red}}_y(q))=n-2,
\end{equation*}
and finally
\begin{equation*}
  \rank(D\! \mbox{\RS}_y(q))=\rank (DR^{\mathrm{red}}_y(q))=n-2.
\end{equation*}
Since rank of $D\! \mbox{\RS}_y(q)$ is maximal, therefore matrix $D\mbox{\RS}_y(q)$ is non-degenerate.
\qed

From Lemmas~\ref{lem:xblk-non-deg}  and~\ref{lem:cblk-y-nondeg} we obtain immediately
\begin{theorem}
\label{thm:collndegRSd2}
Let $d=2$ and   $\qnm$   be a collinear nCC satisfying
 the normalization  $y_{n-1}=0$.
 
 If $x_{n-1} \neq 0$, then $\qnm$ is a non-degenerate solution of \RS.
\end{theorem}

\subsection{Spatial case}
\label{subsec:coll-spatial}

If $d=3$, it may occur that some collinear nCCs are degenerate solutions of \RS.  Observe  that, by (\ref{eq:DRredxy}), the block $D\mbox{\RS}_z(q)$ is obtained from $D\mbox{\RS}_y(q)$ by removing  the row and column corresponding to the $(n-2)$-th body.  Although Lemma~\ref{lem:cblk-y-nondeg} implies that $D\mbox{\RS}_y(q)$ is in isomorphism, the removal
of a row and column with the same index may make the resulting matrix degenerate.  We observed that this happens for $n=5$  for convex combinations
of the following sets of mass parameters 
($m_0=m_1=m_2=0.25, m_3=m_4=0.125$) and ($ m_0=m_1=m_2=m_3=0.166667, \quad m_4=1-(m_0+m_1+m_2+m_3) = 0.333332$), for collinear nCC satisfying the ordering of bodies $x_1 < x_2 < x_0 < x_3 < x_4$. The indexing of bodies
corresponds to the reduced system implemented in our program, with $q_{n-1}$ being computed from the center of mass condition
and $k_1=n-2$ and $k_2=0$ (see system (\ref{eq:cc-red-3D})).
We conjecture the following.

\begin{con}
\label{con:ndeg}
Consider $d=3$. Assume that $\qnm$ is a collinear nCC, such that $y_{n-1}=0$. Then there exist a permutation of bodies 
such that $\qnm$ is a non-degenerate solution of induced \RS.
\end{con}

\section{How we handle degeneracies in the program}
\label{sec:about-prog}

Basic algorithm  is described in Section 7 in \cite{MZ}.  We proceed in two stages: a searching stage and a testing stage. In the searching stage, we cover the set of all possible configurations by cubes and perform successive bisections until, for each cube, we can determine whether it contains a unique zero, contains no zero, or becomes smaller than a prescribed threshold, in which case it is labeled as {\em undecided}. For undecided boxes, we apply additional heuristics to resolve them. This stage may fail, in which case the program yields no conclusion regarding the finiteness of central configurations. In the testing stage, the program identifies the central configurations, since the same configuration may be obtained from different boxes that either overlap or are related by symmetry.

Compared to the programs described in~\cite{MZ,MZ20}, which perform reasonably well in the equal-mass case, the situation with unequal masses is more delicate. In particular, degeneracies that may arise when passing to the reduced system become an issue of significant importance. To address this problem, our program implements two techniques:
\begin{itemize}
\item rotation and/or permutation of the bodies in the configuration, which amounts to choosing a different reduced system;
\item the use of symmetry arguments to reduce the search space, thereby avoiding certain degeneracies.
\end{itemize}

\subsection{Degeneracy conditions}
\subsubsection{Planar case}
 By Theorems~\ref{thm:RS-nondeg-equiv} and~\ref{thm:sol-RS-deg}, when $d=2$, the choice of the reduced system can cause degeneracy only when one of the following conditions holds:
\begin{eqnarray}
x_{n-1}&=&0, \label{eq:xn1=0} \\
  x_{n-1}&=&x_n. \label{eq:xn1xn} 
  \end{eqnarray}

\subsubsection{Spatial case}
For $d=3$ the situation is more involved.  By Theorems~\ref{thm:RS3D-nondeg-equiv} and~\ref{thm:sol-RS3D-deg} 
the degeneracy resulting from the choice of a particular reduced system arises in the following situations:
\begin{itemize}
\item one of conditions  (\ref{eq:xn1xn}) or  (\ref{eq:xn1=0}) is satisfied,
\item the configuration is not collinear and one of the following conditions is satisfied
\begin{eqnarray}
  y_{n-2}=0, \label{eq:yn2=0} \\
   \det \left[\begin{array}{cc}
  (y_{n-2}-y_n) &   (y_{n-1}-y_n)  \\
  (x_{n-2}-x_n)  & (x_{n-1}-x_n)  \\
\end{array}
\right] = 0,  \label{eq:det=0-prog}
\end{eqnarray}
\item the configuration is collinear and (see the discussion in Section~\ref{subsec:coll-spatial})
\begin{equation}
\det D\mbox{\RS}_z(q) = 0.   \label{eq:detDRSz=0}
\end{equation}
\end{itemize}

\subsubsection{Permutations and rotations to avoid the degeneracy in the reduced system }

In the following discussion of degeneracy conditions and the ways to avoid them, we consider conditions expressed in the form $h(q) = 0$, where $q$ denotes a configuration and $h$ is a smooth function. In our program, we work with interval boxes, denoted by $[q]$; consequently, we never have $h([q]) = 0$. Instead, we can check whether $|h([q])| > \delta > 0$. If we cannot guarantee that $|h([q])| > \delta$, we conclude that the condition $h(q) = 0$ may be satisfied for some $q \in [q]$, which may lead to degeneracy. In such a situation, we perform a rotation and/or permutation of the box $[q]$ to obtain a new box $[\bar q]$, such that for each $q \in [q]$ there exists a congruent configuration $\bar q \in [\bar q]$. Let us stress that the permutations which change the reduced system necessarily involve at least one of the bodies $(n-i)$ with $i = 0,1$, and, in the spatial case, additionally the $(n-2)$-th body. Any permutation that involves the $(n-1)$-th body (or, in the spatial case, also the $(n-2)$-th body) must be accompanied by a rotation, which is required to normalize the configuration.
 Each rotation creates a new box $[\tilde q]$ with a larger volume than the original box $[q]$. This happens because of the wrapping effect in interval arithmetic (see \cite{Mo}) and because the rotation angle is itself an interval, with a diameter similar to that of $[q_i]$, where $i$ is the index of the body placed on the $OX$-axis or in the $OXY$-plane. Therefore, we try to avoid rotations whenever possible. 

In the following discussion
 we assume that the box $[q]$ reduces to a single point, $q$. If a procedure avoids the degeneracy condition $s$, the same procedure can be applied to a box with a small diameter, as long as it contains no collisions between bodies.
We consider different orders of the bodies and rotations of the configuration. 

In the sequel by $O_{a}(t)$ we will denote a rotation by an angle $t$ around $Oa$-axis, where $a \in \{x,y,z\}$

 The worst of all degeneracy conditions is condition (\ref{eq:xn1=0}). 
  \begin{center}
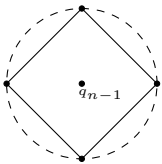

\scalebox{1.0}[1.0]{
\begin{tikzpicture}
\draw[dashed] (0.0, 0.0) circle (1.0cm);

\draw[fill] (0.0, 0.0) circle (0.035cm);
\node[] at (0.25, -0.15) {\tiny $q_{n-1}$};

\draw[fill] (1.0, 0.0) circle (0.035cm);
\draw[fill] (0.0, 1.0) circle (0.035cm);
\draw[fill] (-1.0, 0.0) circle (0.035cm);
\draw[fill] (0.0, -1.0) circle (0.035cm);

\draw[] (1.0, 0.0) -- (0.0, 1.0); 
\draw[] (0.0, 1.0) -- (-1.0, 0.0); 
\draw[] (-1.0, 0.0) -- (0.0, -1.0); 
\draw[] (1.0, 0.0) -- (0.0, -1.0); 

\end{tikzpicture}
}
\captionof{figure}{A continuous family of CCs generated by rotations about the center of mass coinciding with $q_{n-1}$.}. \label{fig:shpere}
\end{center}

If this condition is satisfied, then our box $[q]$ contains in its interior an nCC $q$ (with $x_{n-1}=0$); consequently, it will also contain $O_z(\varphi)q$ for $\varphi$ sufficiently small (note that, due to our normalization, $y_{n-1}=0$ and, in the spatial case, also $z_{n-1}=0)$, see Fig.~\ref{fig:shpere}. This implies that any box containing $q$ will remain undecided if we stick to a fixed reduced system.
To remedy this, we change the body placed on the $OX$-axis. We select the body farthest from the origin and permute the bodies so that this body becomes the $(n-1)$-th one. We then rotate the coordinate system so that this body lies on the $OX$-axis (see Fig~\ref{fig:hard-degen}). In this way, also condition~(\ref{eq:xn1xn}) will  be not satisfied. 

\begin{center}
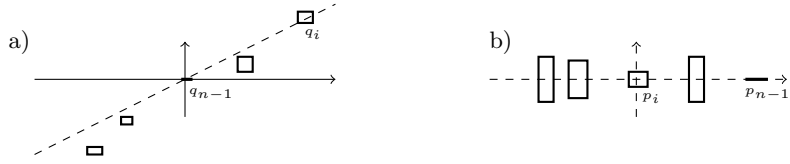

\scalebox{1.0}[1.0]{
\begin{tikzpicture}
\node[] at (-2.2,0.5) {\small a)};
\draw[->] (-2.0, 0.0) -- (2.0, 0.0);
\draw[->] (0.0, -0.5) -- (0.0, 0.5);
\draw[dashed] (-2.0, -1.0) -- (2.0, 1.0);

\draw[thick] (-1.3, -1.0) rectangle (-1.1, -0.9);
\draw[thick] (-0.85, -0.6) rectangle (-0.7, -0.5);
\draw[very thick] (-0.05, 0.0) rectangle (0.1, 0.0); 
\draw[thick] (0.7, 0.1) rectangle (0.9, 0.3);
\draw[thick] (1.5, 0.75) rectangle (1.7, 0.9);

\node[] at (0.35, -0.2) {\tiny $q_{n-1}$};
\node[] at (1.7, 0.6) {\tiny $q_i$};

\node[] at (4.2,0.5) {\small b)};
\draw[dashed,->] (4.05, 0.0) -- (8.0, 0.0);
\draw[dashed,->] (6.0, -0.5) -- (6.0, 0.5);

\draw[thick] (4.7, -0.3) rectangle (4.9, 0.3);
\draw[thick] (5.1, -0.25) rectangle (5.35, 0.25);
\draw[thick] (5.9, -0.1) rectangle (6.15, 0.1); 
\draw[thick] (6.7, -0.3) rectangle (6.9, 0.3);
\draw[very thick] (7.45, 0.0) rectangle (7.75, 0.0); 

\node[] at (6.2, -0.25) {\tiny $p_i$};
\node[] at (7.75, -0.2) {\tiny $p_{n-1}$};
\end{tikzpicture}
}
\captionof{figure}{Example of 
planar 5-body configuration a)  before rotation, $q_{n-1}$ is close to the center of mas b) after rotation, $q_i$ becomes the penultimate $p_{n-1}$ body. After rotation, the interval diameters increase significantly.}\label{fig:hard-degen}
\end{center}

Condition (\ref{eq:xn1xn}) is easily avoided by suitable permutation of bodies.
With this we are dealt with the planar case. 

Condition (\ref{eq:yn2=0}) for configurations that are not nearly collinear is handled as follows (see Fig.~\ref{fig:non-col}). First, we find the body farthest from the $OX$-axis and call it $q_i$. Next, we rotate the configuration around the $OX$-axis to move the $i$-th body onto the $OXY$-plane. Finally, we swap the $i$-th body with the $(n-2)$-th body.

\begin{center}
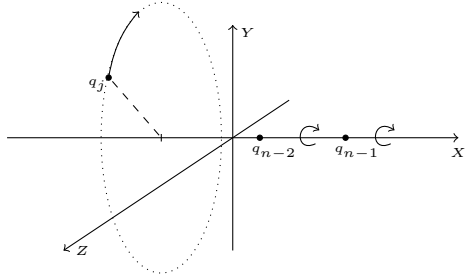

\scalebox{1.0}[1.0]{
\begin{tikzpicture}
\draw[->] (0.0, 1.5) -- (6.0, 1.5);
\node[] at (6.0, 1.3) {\tiny $X$};
\draw[->] (3.0, 0.0) -- (3.0, 3.0);
\node[] at (3.2, 2.9) {\tiny $Y$};
\draw[<-] (0.75, 0.0) -- (3.75, 2.0);
\node[] at (1.0, 0.0) {\tiny $Z$};

\draw[fill] (4.5, 1.5) circle (0.035cm);
\node[] at (4.65, 1.3) {\tiny $q_{n-1}$};

\draw[fill] (3.36, 1.5) circle (0.035cm);
\node[] at (3.55, 1.3) {\tiny $q_{n-2}$};

\draw[fill] (1.35, 2.3) circle (0.035cm);
\node[] at (1.2, 2.2) {\tiny $q_j$};
\draw[dashed] (1.4, 2.25) -- (2.05, 1.5);
\draw[dashed] (2.05, 1.45) -- (2.05, 1.55);

\draw[dotted] (2.05,1.5) ellipse (0.8cm and 1.8cm);
\draw[->] (1.35, 2.3) .. controls (1.43,2.78) and (1.6, 3.0) .. (1.75, 3.18);

\draw[->] (5.1, 1.42) .. controls (4.8, 1.28) and (4.85, 1.8) .. (5.15, 1.6);
\draw[->] (4.1, 1.42) .. controls (3.8, 1.28) and (3.85, 1.8) .. (4.15, 1.6);
\end{tikzpicture}
}
\captionof{figure}{Rotation around the $OX$-axis.}\label{fig:non-col}
\end{center}

Condition (\ref{eq:det=0-prog}) is handled by a suitable permutation of bodies. This condition is equivalent to the following (see Figure~\ref{fig:lin-dep}): the projections of the $(n-1)$-th, $(n-2)$-th, and $n$-th bodies onto the $OXY$-plane are  collinear. 
 Since rotation increases the diameter of the configuration, we prefer, whenever possible, not to change the positions of the ($n-1$)-th and ($n-2$)-th bodies, as doing so would require a rotation. The $n$-th body has no special constraints—it is simply determined by the center-of-mass condition—so we choose a different body to play the role of the computed one.
 Because we assume that conditions~(\ref{eq:xn1=0}) and~(\ref{eq:yn2=0}) are not satisfied, there must exist a body $q_i$ that its projection onto plane $OXY$ does not lie on the line passing through ($x_{n-1}, y_{n-1}$) and ($x_{n-2}, y_{n-2}$); otherwise, the projection of the center of mass would lie on that line, which is possible because   the center of mass  is at the origin of the coordinate system.

\begin{center}
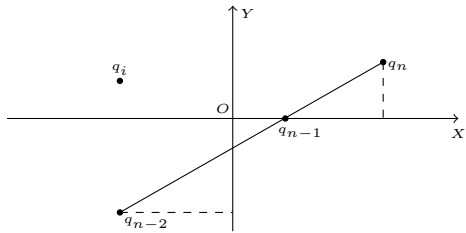

\scalebox{1.0}[1.0]{
\begin{tikzpicture}
\draw[->] (0.0, 1.5) -- (6.0, 1.5);
\node[] at (6.0, 1.3) {\tiny $X$};
\draw[->] (3.0, 0.0) -- (3.0, 3.0);
\node[] at (3.2, 2.9) {\tiny $Y$};
\node[] at (2.87, 1.63) {\tiny $O$};

\draw[] (1.5, 0.25) -- (5.0, 2.25);

\draw[fill] (1.5, 0.25) circle (0.035cm);
\node[] at (1.85, 0.1) {\tiny $q_{n-2}$};
\draw[dashed] (1.5, 0.25) -- (3.0, 0.25);

\draw[fill] (3.7, 1.5) circle (0.035cm);
\node[] at (3.9, 1.3) {\tiny $q_{n-1}$};

\draw[fill] (5.0, 2.25) circle (0.035cm);
\node[] at (5.2, 2.2) {\tiny $q_{n}$};
\draw[dashed] (5.0, 2.25) -- (5.0, 1.5);

\draw[fill] (1.5, 2.0) circle (0.035cm);
\node[] at (1.5, 2.15) {\tiny $q_i$};

\end{tikzpicture}
}
\captionof{figure}{Linearly dependent vectors in condition~(\ref{eq:det=0-prog}); conditions~(\ref{eq:xn1=0}) and~(\ref{eq:yn2=0}) are not satisfied.}\label{fig:lin-dep}
\end{center}

In the case of nearly collinear configuration along $OX$-axis we just try to permute bodies to avoid (\ref{eq:detDRSz=0}), by
changing $(n-2)$-th body, while keeping $(n-1)$-th the same.  Observe that this requires rotation around $OX$-axis to normalize the configuration.

\subsection{Restricting the search space using symmetry arguments}
\label{subsec:RestSearchSet}

In principle, we can fix the indexing of bodies by choosing a reduced system and defining a box in the reduced configuration space ($x_{n-1} = y_{n-1} = z_{n-1} = 0$ and $z_{n-2} = 0$; in the planar case, we ignore the $z$-coordinates), based on a priori bounds as obtained in~\cite{MZ}. We can then run a subdivision algorithm, using the tools described earlier to handle non-degeneracies. In principle, this should work if all nCCs are non-degenerate and Conjecture~\ref{con:ndeg} holds.

 However, this approach turns out to be very inefficient. For the spatial case with $n = 5$ and unequal masses, we were not able to complete a successful run in a reasonable time. 

 The reason can be illustrated as follows.  Consider an equilateral quadrangle with a body at the origin (see Fig.~\ref{fig:shpere}). Assume that the penultimate body $q_{n-1}$ is at the origin and $q_{n-2}$ lies on the $OX$-axis. If we rotate the configuration around the $OZ$-axis, we obtain a circle of nCCs for which $x_{n-1} = y_{n-1} = z_{n-1} = 0$ and $z_{n-2} = 0$. Hence, in the reduced configuration space, this forms a circle of nCCs. In the spatial case, there are also rotations around the $OX$-axis, so our nCCs form part of a two-dimensional set.

 During the algorithm, this continuum of nCCs ends up being covered by a large number of small boxes, which must be rotated and permuted to avoid degeneracies. This process, especially in the spatial case, leads to practical stalling of the program.

The idea is to restrict the search space using symmetry arguments.

The case of equal masses  best illustrates how to avoid the degeneracies described above. In this setting, we may assume that, in a central configuration, the penultimate body on the $OX$-axis is the farthest from the origin, with its $x$-coordinate at least $0.5$ (see Lemma 10 in~\cite{MZ}). We also assume that the ($n-2$)-th body is the farthest from the $OX$-axis, and that there exists an ordering of the $x_i$ for the remaining bodies (see Section 6.1 in~\cite{MZ} or Section 5.1 in \cite{MZ20}, where a slightly different indexing is used). All other configurations can be obtained by permutations of bodies.
With this choice of the search space, we avoid all degeneracies listed above, except for a possible degeneracy corresponding to collinear nCCs in the spatial case. For $n = 5,6$, this situation does not occur, as the successful runs of our program reported   in \cite{MZ20} show.

In the unequal mass case, instead of a single run, we perform $n$ runs, each time placing a different body on the positive $OX$-axis and assuming that this body is the farthest from the origin. Under this assumption, we have $x_{n-1} \ge 0.5$ and $x_{n-1} > x_n$, that is, the negations of conditions~\ref{eq:xn1xn}) and~(\ref{eq:xn1=0}) hold. In the planar case, these assumptions eliminate all possible degeneracies of the reduced systems. In the spatial case, however, we still need to handle the remaining degeneracies.

For unequal masses in the planar 5-body case we run the program in five separate runs, each time placing a different body on the $OX$-axis as the farthest one. However, due to the absence of additional constraints present in the equal-mass case, for certain difficult mass distributions the runtime is not merely five times longer (as one might expect from performing five runs); instead, it increases by a factor of approximately 114.

When some of the masses are equal, we can reduce the number of runs. In the extreme case of equal masses, a single run is sufficient.

\section{Some of exceptional cases for five bodies}\label{sec:diff-masses}
In this section, we consider several challenging mass distributions that arise in the work of~\cite{AK} in two dimensions and in the work of~\cite{HJ} in three dimensions.

\subsection{One exceptional and two 'difficult' cases in 2D}\label{sec:except-case-2D}

In~\cite{AK}, exceptional cases of mass parameters are identified, corresponding to the diagrams in Figure 11 of that work. For these cases, the finiteness of the number of equivalence classes of central configurations has not been established.

In the present work, we investigate one of these exceptional cases. Although this cases are defined by $m_1=m_2$ and $m_3=m_4$, our analysis is restricted to several specific discrete choices of mass parameters. Consequently, we do not claim to resolve the finiteness question for this exceptional case in general.  We also investigate  two other cases  of polynomial relation between masses defined by a single equation, hence not exceptional
 in the sense of Albouy-Kaloshin, but they are derived from some difficult diagrams in Figure 11 of  \cite{AK}.  For these cases we obtain exceptional masses, when some other polynomial equation (there are many of them) is added, as  explained in the proof of Theorem 6 on page 582 of \cite{AK}.

When presenting our results for different mass values, we classify planar CCs as collinear, concave, or convex (pentagonal). Collinear solutions, which arise from Moulton’s Theorem and Conley’s Theorem, are treated separately from the other cases; this allows our program to classify them unambiguously.

For concave CCs, we initially attempted a finer classification based on the number of bodies lying in the interior of the convex hull, distinguishing between triangular and quadrilateral hulls. However, we observed that as the mass parameters vary, some CCs evolve continuously from one such class to another without undergoing any bifurcation. Moreover,  even the distinction between concave and convex CCs may not be sharp for the $n$-body problem when $n \geq 5$. While for $n=4$ the comprehensive results about existence concave central configurations has been given first in \cite{H02} and  it is known (see~\cite{MB,X04}) that any convex planar central configuration for the Newtonian four-body problem must be strictly convex, an analogous result does not hold for $n=5$. Indeed, \cite{CH} provides explicit examples of central configurations that are convex but not strictly convex. Notably, the examples in~\cite{CH} also belong to the exceptional cases considered in this section.

The cases we investigated  do not pose serious difficulties for our program. While there are various ways to prove that a given configuration contains no solution, only the Newton–Krawczyk method allows us to rigorously certify the existence of a solution. Consequently, the only genuinely challenging situations arise at bifurcation points, where the Jacobian matrix fails to be an isomorphism and the Newton–Krawczyk method is not applicable. Our particular choices of mass parameters avoid these bifurcation loci.

Below we summarize the results obtained for individual mass configurations for three cases.   Each of these cases 
 is  a manifold, however we just treat a sample of points. The varying number of solutions suggests the presence of bifurcation points somewhere between the sampled mass values—these would likely correspond to the difficult cases for our algorithm. For all configurations, the number of collinear solutions and pentagon configurations remains constant, whereas the number of concave configurations changes with the masses.  Moulton's theorem gives $60$ collinear CCs and we have $24=4!$ pentagons corresponding to a different cyclic order of bodies.  We obtained numbers of central configurations between 294 and 450, as suggested by results of  Simo \cite{Si03} (see also \cite{A15}). 

\begin{enumerate}
\item\label{enum:1} Case $m_1 = m_2$, $m_3 = m_4$; this is equation (32) (case 8.6)  in \cite{AK}. This is an exceptional case of mass parameters in \cite{AK}.

For the sample masses, the program finds the following number of solutions of a given type (shown in Table~\ref{table:0}). Note that the equal-mass case is a special case, which is also included in the table.

\begin{table}[H]
\centering
\begin{tabular}{p{3.6cm}|p{2cm}|p{2cm}|p{2cm}|p{2cm}}
 & concave & collinear & pentagons & \textbf{total}\\
\hline
$\begin{array}{lcl}m_1 = m_2 & = &\\ m_3 = m_4 & = & 0.2\end{array}$ & 270 &  60 & 24 & \textbf{354}\\
\hline
$\begin{array}{lcl}m_1 = m_2 & = & 0.21\\ m_3 = m_4 & = & 0.19\end{array}$ & 270 & 60 & 24 & \textbf{354}\\
\hline
$\begin{array}{lcl}m_1 = m_2 & = & 0.22\\ m_3 = m_4 & = & 0.18\end{array}$ & 246 & 60 & 24 & \textbf{330}\\
\hline
$\begin{array}{lcl}m_1 = m_2 & = & 0.3\\ m_3 = m_4 & = & 0.1\end{array}$ & 218 & 60 & 24 & \textbf{302}\\
\hline
$\begin{array}{lcl}m_1 = m_2 & = & 0.35\\ m_3 = m_4 & = & 0.05\end{array}$ & 226 & 60 & 24 & \textbf{310}\\
\hline
$\begin{array}{lcl}m_1 = m_2 & = & 0.39\\ m_3 = m_4 & = & 0.01\end{array}$ &  242 & 60 & 24 & \textbf{326}\\
\hline
$\begin{array}{lcl}m_1 = m_2 & = & 0.2499\\ m_3 = m_4 & = & 0.2499 \end{array}$ &  318 & 60 & 24 & \textbf{402}\\
\hline
$\begin{array}{lcl}m_1 = m_2 & = & 0.24999\\ m_3 = m_4 & = & 0.24999 \end{array}$ &  366 & 60 & 24 & \textbf{450}\\
\end{tabular}
\caption{Summary of the number of distinct solutions for five bodies under the equal-mass-pair criterion  \textbf{($m_1 = m_2$, $m_3 = m_4$)} on the plane.}\label{table:0}
\end{table}

\item Case $m_1 m_3 = m_2 m_4$; this is equation (29) (case 8.3) in \cite{AK}.    For the sample masses, the program finds the following number of solutions of a given type see Table~\ref{table:1}. Notice that
all cases presented in Table~\ref{table:0} also fulfill the corresponding condition $m_1 m_3 = m_2 m_4$; consequently, Table~\ref{table:0} provides alternative sample results.  Observe that 
equation $m_1 m_3 = m_2 m_4$ defines codimension one subvariety in the mass space, hence not all points on it are exceptional.  To obtain an exceptional case one needs to add another polynomial relation
as explained  in the proof of Theorem 6 on page 582 of \cite{AK}.

\begin{table}[H]
\centering
\begin{tabular}{p{2.5cm}|p{2cm}|p{2cm}|p{2cm}|p{2cm}}
 & concave & collinear & pentagon & \textbf{total}\\
\hline
$\begin{array}{lcl}
m_1 & = & 0.22\\
m_2 & = & 0.18\\
m_3 & = & 0.11\\
m_4 & = & 0.36
\end{array}$ & 210 & 60 & 24 & \textbf{294} \\
\hline
$\begin{array}{lcl}
m_1 & = & 0.22\\
m_2 & = & 0.18\\
m_3 & = & 0.22\\
m_4 & = & 0.18
\end{array}$ &  246 & 60 & 24 & \textbf{330}\\
\hline
$\begin{array}{lcl}
m_1 & = & 0.22\\
m_2 & = & 0.18\\
m_3 & = & 0.33\\
m_4 & = & 0.12
\end{array}$ & 210 &  60 & 24 & \textbf{294}\\
\hline
$\begin{array}{lcl}
m_1 & = & 0.22\\
m_2 & = & 0.18\\
m_3 & = & 0.44\\
m_4 & = & 0.09
\end{array}$ & 210 &  60 & 24 & \textbf{294}\\
\end{tabular}
\caption{Summary of the number of distinct solutions for five bodies under the product criterion \textbf{($m_1 m_3 = m_2 m_4$)} on the plane.}\label{table:1}
\end{table}

\item\label{item:num1} Case $\frac{1}{\sqrt{m_3}} = \frac{1}{\sqrt{m_1}} + \frac{1}{\sqrt{m_2}}$; this is equation (23) (case 8.1) in \cite{AK}.

In this case, the mass of the third body calculated from the root equation is an interval; in the Table~\ref{table:2}, we only give the first four digits, which are the same for the left and right ends of the range.

Observe that 
equation $\frac{1}{\sqrt{m_3}} = \frac{1}{\sqrt{m_1}} + \frac{1}{\sqrt{m_2}}$ defines a codimension one subvariety in the mass space, hence not all points on it are exceptional.  To obtain an exceptional case one needs to add another polynomial relation
as explained  in the proof of Theorem 6 on page 582 of \cite{AK}.

\begin{table}[ht]
\centering
\begin{tabular}{p{3cm}|p{2cm}|p{2cm}|p{2cm}|p{2cm}}
 & concave & collinear & pentagon & \textbf{total} \\
\hline
$\begin{array}{lcl}
m_1 & = & 0.3\\
m_2 & = & 0.3\\
m_3 & = & 0.075\\
m_4 & = & 0.21
\end{array}$ & 218 & 60 & 24 & \textbf{302} \\
\hline
$\begin{array}{lcl}
m_1 & = & 0.3\\
m_2 & = & 0.3\\
m_3 & = & 0.075\\
m_4 & = & 0.22
\end{array}$ & 214 & 60  & 24 & \textbf{298}\\
\hline
$\begin{array}{lcl}
m_1 & = & 0.3\\
m_2 & = & 0.25\\
m_3 & = & 0.0683\\
m_4 & = & 0.21
\end{array}$ & 222  & 60 & 24 & \textbf{306} \\
\hline
$\begin{array}{lcl}
m_1 & = & 0.3\\
m_2 & = & 0.25\\
m_3 & = & 0.0683\\
m_4 & = & 0.22
\end{array}$ & 222 & 60 & 24 &  \textbf{306}\\
\hline
$\begin{array}{lcl}
m_1 & = & 0.3\\
m_2 & = & 0.2\\
m_3 & = & 0.0606\\
m_4 & = & 0.21
\end{array}$ & 222 & 60  & 24 & \textbf{306} \\
\hline
$\begin{array}{lcl}
m_1 & = & 0.3\\
m_2 & = & 0.2\\
m_3 & = & 0.0606\\
m_4 & = & 0.22
\end{array}$ & 222 & 60  & 24 & \textbf{306}\\
\end{tabular}
\caption{Summary of the number of distinct solutions for  five bodies under the square-root criterion ($\frac{1}{\sqrt{m_3}} = \frac{1}{\sqrt{m_1}} + \frac{1}{\sqrt{m_2}}$) on the plane.}\label{table:2}
\end{table}

\end{enumerate}

To conclude, we discuss in more detail the examples presented in \cite{CH}. The authors consider the planar five-body problem with masses that are not normalized to satisfy $\sum_{i=1}^n m_i = 1$. In their examples, the masses satisfy the relations defining our exceptional cases 1 and 2 above. They establish the existence of central configurations that are convex but not strictly convex. For such configurations, our program would be unable to determine whether the configuration is convex.

We now briefly summarize their results. For masses
$$
m_1 = m_2 = m_5 = 1, \qquad m_3 = m_4 = \mu \approx 11.23156072828415553841745,
$$
the authors prove the existence of a central configuration that is a local minimum of the normalized potential $\sqrt{I}U$.

In addition, other numerical solutions described in \cite{CH} with
$$
m_1 = m_2 = 1, \qquad m_3 = m_4 = \mu, \qquad m_5 = \nu,
$$
suggest the presence of a one-parameter family of solutions in which both $\mu$ and $\nu$ increase simultaneously. Two representative examples are ($\nu = 10^{-4}, \mu = 2.758$), which does not correspond to a local minimum, and ($\nu = 10^{4}, \mu = 81952.332$), which does. For $\nu \approx 0.5180855751$, the solution becomes degenerate. Note that in this degenerate case, our program would simply fail to establish finiteness.

\subsection{Exceptional cases in 3D}
In \cite{HJ}, Hampton and Jensen establish the finiteness of spatial, non-planar central configurations for the five-body problem, with the exception of explicitly described special cases of mass values. Their result generalizes an earlier generic finiteness theorem of Moeckel \cite{M01}. The exceptional cases are listed in Table 1 of \cite{HJ}, where each case corresponds to a row in the table and is identified by its first entry, referred to in \cite{HJ} as the {\em ray index}.

In the present work, we investigate two specific representatives of two exceptional cases corresponding to the first and third rows of Table 1 in \cite{HJ}.

For the first row of Table 1 in \cite{HJ}, the ray index is [59], and the exceptional polynomial is
\begin{equation}\label{eq:r59-1}
m_1 m_2 - m_3 m_4 - m_3 m_5 = 0.
\end{equation}

For the third row of Table 1 in \cite{HJ}, the ray index is [59,72], and there are two exceptional polynomials,
\begin{eqnarray}
m_3 - m_4 - m_5 &=& 0, \label{eq:r59-2}\\
m_4^2 + 2 m_4 m_5 + m_5^2 - m_1 m_2 &=& 0. \label{eq:r59-3}
\end{eqnarray}

We consider two sets of mass parameters:
$$
(m_1,m_2,m_3,m_4,m_5)=(2,1,1,1,1),
$$
which satisfies \eqref{eq:r59-1}, and
$$
(m_1,m_2,m_3,m_4,m_5)=(2,2,2,1,1),
$$
which satisfies all relations \eqref{eq:r59-1}–\eqref{eq:r59-3}. In our computations, the masses are normalized so that $\sum m_i=1$. Note that these mass choices also correspond to exceptional cases 1 and 2 for the planar system discussed in Section~\ref{sec:except-case-2D}.

Computing directly all central configurations  for these mass distributions is  computationally very expensive. Therefore, for each of the above mass distributions, we restrict ourselves to two representative runs (instead of five needed in general case): one in which a heavy body is placed on the $OX$-axis and one in which a light body is placed on the $OX$-axis. All remaining runs are equivalent by symmetry and would yield identical results. Although this procedure does not produce all possible CCs directly, all others can be recovered via symmetry arguments. Consequently, we obtain finiteness of the number of central configurations. While it would be possible, with additional effort, to compute the exact number of equivalence classes, in this work we restrict ourselves to providing upper bounds.

For the mass distribution ($2,1,1,1,1$), the program produces 106\ 210 solutions in the run with a light body placed on the $OX$-axis, while no solutions are found when the heavy body is placed on the $OX$-axis. After identification, these solutions correspond to 94 distinct central configurations, including 18 collinear, 42 non-collinear planar, and 34 spatial non-planar configurations. In this case, the upper bound on the number of central configurations is $4 \cdot 94$, rather than $5 \cdot 94$, since one of the runs yields no solutions. These data are summarized in Table~\ref{tab:CC3D} as the first column.

For the mass distribution ($2,2,2,1,1$), the program finds 52 038 solutions when a light body is placed on the $OX$-axis and 412 solutions when a heavy body is placed on the $OX$-axis. After identification these reduce to 157 distinct central configurations, consisting of 30 collinear, 77 non-collinear planar, and 50 spatial non-planar configurations. Since 157 provides an upper bound for configurations in which a given body is the farthest from the origin and is placed on the $OX$-axis, the overall upper bound on the number of central configurations is $785 = 5 \cdot 157$. These data are summarized in Table~\ref{tab:CC3D} as the second column.

\begin{table}[H]
\centering
\begin{tabular}{lcc}
\hline
 & $2,1,1,1,1$ & $2,2,2,1,1$ \\
\hline
Solutions (light b. on $OX$) & 106{,}210 & 52{,}038 \\
Solutions (heavy b. on $OX$) & 0 & 412 \\
Distinct CCs & 94 & 157 \\
Collinear CCs & 18 & 30 \\
Planar (non-collinear) CCs & 42 & 77 \\
Spatial (non-planar) CCs & 34 & 50 \\
\hline
\textbf{upper bound on \#CCs} & \textbf{376} & \textbf{785}\\
\hline
\end{tabular}
\caption{Summary of solutions and distinct central configurations (CCs) for the two mass distributions for five bodies in the spatial case. Notice that the numbers of solutions are obtained from only two runs; to obtain the exact numbers we should use an upper bound with a symmetry arguments.}
\label{tab:CC3D}
\end{table}

Let us note that the case of equal masses is also exceptional in the spatial setting, since it satisfies relation~(\ref{eq:r59-1}). This case was treated in \cite{MZ20}. As in the planar case, a single computational run is sufficient, with all other central configurations obtained via symmetry arguments. The number of distinct equivalence classes of central configurations is equal to $307$ (see Table~\ref{tab:eq-CC3D}).

Observe that this number is smaller than the number of non-equivalent planar central configurations, which is 354 (see the first row of Table~\ref{table:0}). This difference arises because, in three dimensions, certain planar central configurations (see Fig.~\ref{fig:square}) that are distinct under the action of $SO(2)$ become equivalent when the action of $SO(3)$ is taken into account.

 \begin{center}
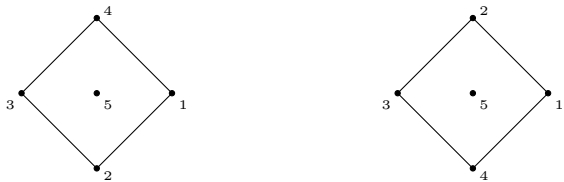

\scalebox{1.0}[1.0]{
\begin{tikzpicture}
\draw[fill] (0.0, 0.0) circle (0.035cm);
\node[] at (0.15, -0.15) {\tiny $5$};

\draw[fill] (1.0, 0.0) circle (0.035cm);
\node[] at (1.15, -0.15) {\tiny $1$};

\draw[fill] (0.0, 1.0) circle (0.035cm);
\node[] at (0.15, 1.1) {\tiny $4$};

\draw[fill] (-1.0, 0.0) circle (0.035cm);
\node[] at (-1.15, -0.15) {\tiny $3$};

\draw[fill] (0.0, -1.0) circle (0.035cm);
\node[] at (0.15, -1.1) {\tiny $2$};

\draw[] (1.0, 0.0) -- (0.0, 1.0); 
\draw[] (0.0, 1.0) -- (-1.0, 0.0); 
\draw[] (-1.0, 0.0) -- (0.0, -1.0); 
\draw[] (1.0, 0.0) -- (0.0, -1.0); 

\draw[fill] (5.0, 0.0) circle (0.035cm);
\node[] at (5.15, -0.15) {\tiny $5$};

\draw[fill] (6.0, 0.0) circle (0.035cm);
\node[] at (6.15, -0.15) {\tiny $1$};

\draw[fill] (5.0, 1.0) circle (0.035cm);
\node[] at (5.15, 1.1) {\tiny $2$};

\draw[fill] (4.0, 0.0) circle (0.035cm);
\node[] at (3.85, -0.15) {\tiny $3$};

\draw[fill] (5.0, -1.0) circle (0.035cm);
\node[] at (5.15, -1.1) {\tiny $4$};

\draw[] (6.0, 0.0) -- (5.0, 1.0); 
\draw[] (5.0, 1.0) -- (4.0, 0.0); 
\draw[] (4.0, 0.0) -- (5.0, -1.0); 
\draw[] (6.0, 0.0) -- (5.0, -1.0); 

\end{tikzpicture}
}
\captionof{figure}{An example of configurations that are identified as one in the spatial case, while they are two different configurations in the planar case.} \label{fig:square}
\end{center}

We also note that the number of central configurations reported here differs from the values previously given in \cite{MZ20}, specifically those listed in Table~4 under the column labeled $\cong(n)$, which is intended to represent the number of equivalence classes under the action of the group $SO(3)$. The values $\cong(4)$ and $\cong(5)$ reported there are incorrect, and it is likely that $\cong(6)$ is also incorrect. The correct values are $\cong(4)=33$ and $\cong(5)=307$. In contrast, the entries $\iso(4)$ and $\iso(5)$ reported in \cite{MZ20} are correct.

\begin{table}[H]
\centering
\begin{tabular}{lc}
\hline
 & equal masses  \\
\hline
Distinct CCs & 307   \\
Collinear CCs & 60  \\
Planar (non-collinear) CCs & 147  \\
Spatial (non-planar) CCs & 100  \\
\hline
\end{tabular}
\caption{Summary of 5-bodies spatial central configurations with equal masses.}
\label{tab:eq-CC3D}
\end{table}

\section{Implementation notes}

Supplementary materials have been attached to the paper \cite{MZ26}, including the program source code, input files, and report files generated during the program execution. The program is written in C++ and has been tested under macOS 15.6.1 using Apple LLVM with clang
version 17.0.0 (clang-1700.6.3.2). The program uses the CAPD library (as of November 2023; see also
the sourceforge download zone).  The execution times of individual computations are considerable. For instance, for equal mass pairs $m_0 = m_1 = 0.21$ and $m_2 = m_3 = 0.19$, the program runs for  4513,95 s. In the spatial case, the program takes even longer to compute. Note that in the program the bodies are numbered starting from 0.

\appendix
\section{One lemma about elimination of variables and rank of the equation}\label{sssec:lm-rank}

Assume that we have variables $x \in \mathbb{R}^{d_1}$ and $y \in \mathbb{R}^{d_2}$ and a system of equations
\begin{subequations}
\begin{align}
  F_1(x,y)&=0, \label{eq:F1}  \\
  F_2(x,y)&=0, \label{eq:F2}
\end{align}\label{eq:F1and2}
\end{subequations}
\noindent
where $F_1: \mathbb{R}^{d_1} \times \mathbb{R}^{d_2} \to \mathbb{R}^{m}$ and $F_2: \mathbb{R}^{d_1} \times \mathbb{R}^{d_2} \to \mathbb{R}^{d_2}$ are $C^1$ functions. Let us stress that the number of equations in $F_2$ agrees with the dimension of $y$.
Let us denote $$F(x,y)=(F_1(x,y),F_2(x,y)).$$
If we can locally eliminate $y$ from equation $F_2(x,y)=0$ i.e. solve for $y$ for a given $x$, to obtain function $y(x)$, then  a reduced system of equations is defined by
\begin{equation}
F_\mathrm{red}(x)=F_1(x,y(x)).
\end{equation}

\begin{lemma}
\label{lem:rank-after-elim}
  Assume that $(x_0,y_0)$ satisfies (\ref{eq:F1and2}) and $\frac{\partial F_2}{\partial y}(x_0,y_0)$ is an isomorphism. Then
  \begin{equation}
     \rank\left(D\! F_\mathrm{red}(x_0)\right) = \rank\left(D\! F(x_0,y_0)\right)-d_2.
  \end{equation}
\end{lemma}
\proof
Since $\frac{\partial F_2}{\partial y}(x_0,y_0)$ is an isomorphism,  we can apply the implicit function theorem to locally eliminate variable $y$ by solving $F_2(x,y)=0$ in the neighborhood
of $(x_0,y_0)$. We obtain $y(x) \in C^1$, such that
\begin{equation}
  \frac{\partial y}{\partial x}(x)= - \left(\frac{\partial F_2}{\partial y}(x,y(x))\right)^{-1}\left(\frac{\partial F_2}{\partial x}(x,y(x)) \right).  \label{eq:derIFT}
\end{equation}
Obviously $y_0=y(x_0)$.
Observe that from (\ref{eq:derIFT}) we obtain
\begin{eqnarray*}
  D\! F_\mathrm{red}(x_0) & =& \frac{\partial F_\mathrm{red}}{\partial x}(x_0)=  \frac{\partial F_1}{\partial x}(x_0,y_0)+   \frac{\partial F_1}{\partial y}(x_0,y_0)\frac{\partial y}{\partial x}(x_0)\\
  & = &
  \frac{\partial F_1}{\partial x}(x_0,y_0) -   \frac{\partial F_1}{\partial y}(x_0,y_0) \left(\frac{\partial F_2}{\partial y}(x_0,y_0)\right)^{-1}\frac{\partial F_2}{\partial x}(x_0,y_0).
\end{eqnarray*}

Let us write $D\! F(x_0,y_0) $ as block matrix
\begin{equation*}
  D\! F(x_0,y_0) =
        \begin{bmatrix}
          \frac{\partial F_1}{\partial x}(x_0,y_0) & \frac{\partial F_1}{\partial y}(x_0,y_0)  \\
           \frac{\partial F_2}{\partial x}(x_0,y_0) & \frac{\partial F_2}{\partial y}(x_0,y_0)
        \end{bmatrix}
 =:
        \begin{bmatrix}
          A_{11} & A_{12}  \\
          A_{21} & A_{22}
        \end{bmatrix}.
\end{equation*}
From the assumption it follows that $A_{22}$ is invertible. Observe that the following matrix has the same rank as $ D\! F(x_0,y_0)$
\begin{eqnarray*}
  B:= D\! F(x_0,y_0) \cdot \begin{bmatrix}
                                  I_{d_1} & 0 \\
                                  C_{21} & A_{22}^{-1}
                               \end{bmatrix}   =
                               \begin{bmatrix}
                                  A_{11} + A_{12}C_{21} & A_{12}A_{22}^{-1} \\
                                  A_{21} + A_{22}C_{21} & I_{d_2}
                                     \end{bmatrix}
\end{eqnarray*}
where $I_{d}: \mathbb{R}^{d}\to \mathbb{R}^{d}$ is the identity and $C_{21}: \mathbb{R}^{d_1} \to \mathbb{R}^{d_2}$ is a linear map to be determined below.
Indeed the rank of $B$ is the same as rank of $D\! F(x_0,y_0)$  because the matrix representing the right factor in the above multiplication is an isomorphism.

Observe that if we set
\begin{equation*}
   C_{21}=-A_{22}^{-1}A_{21},
\end{equation*}
then
\begin{equation*}
  B= \begin{bmatrix}
             A_{11} - A_{12}A_{22}^{-1}A_{21} & A_{12}A_{22}^{-1} \\
              0 & I_{d_2}
           \end{bmatrix}.
\end{equation*}
It is immediate that
\begin{eqnarray*}
  \rank (B) & = & \rank \left( \begin{bmatrix}
             A_{11} - A_{12}A_{22}^{-1}A_{21} & 0 \\
              0 & I_{d_2}
           \end{bmatrix} \right) \\[1ex]
           & = & \rank(  A_{11} - A_{12}A_{22}^{-1}A_{21}) + d_2.
\end{eqnarray*}
Let us now rewrite $ D\! F_\mathrm{red}(x_0)$ in terms of $A_{ij}$. We have
\begin{eqnarray*}
   D\! F_\mathrm{red}(x_0)= A_{11} -  A_{12} A_{22}^{-1}A_{21}.
\end{eqnarray*}

\qed

\end{document}